\documentclass[submission,copyright,creativecommons]{eptcs}
\providecommand{\event}{ACT 2026}

\usepackage{iftex}

\ifpdf
  \usepackage{underscore}
  \usepackage[T1]{fontenc}
\else
  \usepackage{breakurl}   
\fi

\usepackage[utf8]{inputenc}
\usepackage{url}
\usepackage[numbers]{natbib}
\usepackage{graphicx}
\usepackage{tikz}
\usetikzlibrary{arrows.meta, positioning}
 \usepackage{booktabs} 

\title{Mapping Category Theory:\\ Growth, Community Structure, and Knowledge Flow}

\author{Eugenio Llanos
\institute{Universidad de Pamplona\\
Pamplona, Colombia}
\email{eugenio.llanosbal@unipamplona.edu.co}
\and
Julián Angarita
\institute{Grupo de Investigación Estudios Sociales del Conocimiento\\
Universidad de Antioquia\\
Medellín, Colombia}
\email{julian.angarita@udea.edu.co}
\and
Wilmer Leal
\institute{Department of Mechanical \& Aerospace Engineering\\
University of Florida\\
Gainesville, USA}
\email{wleal@ufl.edu}
\and
Valeria de Paiva
\institute{Topos Institute \\
Berkeley, USA}
\email{valeria@topos.institute}
}
\def\titlerunning{Mapping Category Theory}
\def\authorrunning{E. Llanos, J. Angarita, W. Leal, \& V. de Paiva}
\begin{document}
\maketitle

\begin{abstract}
Category theory has evolved from a unifying language for mathematics into a conceptual framework that connects diverse mathematical and scientific disciplines. In this paper, we investigate its historical development and patterns of knowledge flow through a scientometric analysis of 30{,}854 category theory publications indexed in zbMATH and authored by 13,392 researchers. We analyze long-term trends in publication activity, topical diversification, collaboration, citation patterns, and interactions between Mathematical Subject Classification (MSC) topics. The results reveal sustained growth, increasing thematic diversification, and a layered research community characterized by a highly productive core and a small group of interdisciplinary researchers who connect category theory with numerous mathematical and applied domains. Topic-level analyses identify theoretical computer science as the principal applied interface of category theory, followed by quantum theory, with additional but smaller interactions involving systems theory, artificial intelligence, biology, economics, and network science. The results further show that these disciplines connect to category theory through distinct categorical interfaces rather than a single generic topic: theoretical computer science is associated primarily with foundational topics such as \emph{Categories and theories} and \emph{General theory of categories and functors}, whereas quantum theory is centered on \emph{Categorical structures}, \emph{Monoidal categories and operads}, and \emph{Homological algebra in category theory}. Citation patterns further suggest that the conceptual core of category theory is continuously absorbing and abstracting ideas originating in surrounding fields while providing general methods that can propagate across disciplines.
\end{abstract}

\tableofcontents

\section{Introduction}

If we think structurally, the history of category theory is less a linear sequence of results and more a dynamic network of ideas, with phases of consolidation, expansion, and reorganization.

Category theory emerged in the mid-twentieth century not as a foundational program, but as a unifying language designed to facilitate the transfer of constructions across mathematical domains. Introduced by Eilenberg and Mac~Lane in 1945~\cite{eilenberg1945}, the notions of categories, functors, and natural transformations arose from problems in algebraic topology, particularly the need to compare homology theories defined in different contexts. At this stage, category theory functioned primarily as a \emph{translation layer}, enabling results and constructions to move between algebra, topology, and geometry.

A major shift occurred in the 1950s and 1960s with the work of Grothendieck, who placed categorical notions—such as abelian categories, derived functors, sites, and toposes—at the center of algebraic geometry. Category theory thereby evolved from a descriptive language into a structural framework capable of organizing entire research programs. Mathematical objects were increasingly studied via their embedding into abstract categorical environments, from which general results could be transported back to concrete settings.

At the same time, category theory underwent a process of internal development. Core constructions such as limits, adjunctions, and monads~\cite{MacLane1971,Awodey2018} were themselves characterized categorically, leading to the systematic study of enriched, internal, and higher-dimensional categories. In this sense, category theory became a meta-theory: a language for describing not only mathematical objects, but also the formal systems in which those objects are studied.

Over time, category theory expanded into a constellation of interconnected subfields. It provided semantic frameworks for logic through categorical logic and topos theory~\cite{Lawvere1963,LawvereTierney1970,Johnstone1977,LambekScott1986,depaiva2013}, and for computation through the categorical semantics of programming languages~\cite{Moggi1991}. Structures such as cartesian closed categories, monads, and monoidal categories became central tools for modeling computation, compositionality, and resource sensitivity~\cite{Girard1987,seely1989}. More recently, higher category theory and homotopical methods have reshaped areas such as algebraic geometry and type theory, while applied category theory has extended categorical methods to domains including networks, databases, and systems theory~\cite{baez2022applied}.

In parallel with these developments, category theory has increasingly interacted with digital mathematics. Large-scale formal libraries,  proof assistants, and mathematical knowledge databases have transformed how mathematical knowledge is produced and circulated. In this context, category theory functions not only as a mathematical discipline but also as a conceptual framework for organizing mathematical data and relationships. 

Viewed through this lens, the development of category theory can be understood as a sequence of reorganizations in the structure and transmission of mathematical knowledge. Its repeated emergence in contexts of increasing conceptual complexity suggests that it plays a distinctive role as an architecture for managing abstraction in modern mathematics and the sciences.

Yet, despite this rich historical picture, we still lack a systematic, large-scale, and data driven account of the structures through which categorical knowledge has developed and circulated.  How category theory has grown and diversified as research field, how its scientific community is organized, and how categorical knowledge interacts with and circulates across disciplinary boundaries have not been systematically mapped.  Understanding these patterns of development and knowledge flow is the problem we address in this article.

\subsection{History of Ideas and Data Availability}

Recent work on the history of ideas increasingly relies on large-scale textual and bibliographic data~\cite{hall2008,computational-history-2013,leal2022exploration}, but the availability of such data varies significantly across scientific communities. In computational linguistics, for example, the ACL Anthology provides a near-complete, open-access record of the field’s literature, enabling detailed historical and conceptual analyses~\cite{gildea2018}. In contrast, chemistry and materials science have developed knowledge infrastructures ranging from open-access but still far from comprehensive resources, such as the \emph{Open Reaction Database}~\cite{ORD-JACS-2021}, to comprehensive but proprietary databases, such as \emph{Reaxys} and \emph{CAS SciFinder}, which provide integrated environments for searching, analyzing, and working with chemical data. In either case, these resources are expert-curated and supported by mature ontological infrastructures, enabling large-scale analyses of scientific knowledge~\cite{Llanos12660,Weston2019,chemspace-periodic-system-2022}.

Category theory occupies a different position. As a relatively young and highly abstract field~\cite{eilenberg1945}, it lacks both the comprehensive open corpora of computational linguistics and the extensive ontological infrastructure of older, industrially driven disciplines. Its publication landscape is more fragmented, and its conceptual structure is expressed primarily through formal and textual means. This makes it both a challenging and a revealing domain for large-scale analysis.

\subsection{Why Category Theory?}

Despite these limitations, the history described above makes category theory particularly well suited for studying large-scale patterns in the development and transmission of mathematical knowledge. If its role in connecting mathematical domains has shaped the evolution of the field, these interactions should have left observable traces in its bibliographic record. This makes it possible to ask, at scale, how category theory has grown and diversified since its origins, how its research community has developed, and how its connections with other mathematical and applied domains have evolved over time.

There is also a second, forward-looking reason for studying category theory. The field is characterized by a relatively small set of core concepts---such as functors, natural transformations, adjunctions, and universal properties~\cite{MacLane1971}---that recur systematically across subfields and over time. Conceptual innovation often takes the form of controlled generalization and definitional refinement, making structural patterns comparatively explicit and traceable. Moreover, category theory emphasizes relationships over objects, making it naturally compatible with graph-based and relational representations. Definitions frequently encode dependencies through universal properties, equivalences, and compositional structures, all of which can, in principle, be recovered from textual and symbolic data. These features make the large-scale analysis pursued here a natural starting point for what we term \emph{network mathematics}: a broader program aimed at combining external, data-driven perspectives with models of the internal semantic and relational organization of mathematical knowledge, in order to study how mathematical concepts are structured, connected, transmitted, and transformed over time.

\subsection{Bibliographic Traces of a Mathematical Community}
Scientific activity leaves bibliographic traces. Publications, authorship, subject classifications, collaborations, and citations record different aspects of how research communities develop, organize, and communicate knowledge over time. Although individually these records provide only partial information, collectively they form a structured historical footprint of scientific activity. Scientometric methods provide a way to analyze this footprint at scale, revealing patterns that are difficult to perceive from individual publications or historical accounts alone~\cite{behrens2011,visualizingKnowledgeDomains2003,price1975science,Llanos12660,Boyack2005}.

In this work, we use these data not as indicators for evaluating researchers or research quality~\cite{caveats-for-the-use-of-citation,Hicks2015}, but as empirical evidence for studying the structure and evolution of the category theory community. Their resolution is necessarily limited: bibliographic relations can reveal patterns of growth, collaboration, disciplinary organization, and interaction, but they do not directly capture the mathematical concepts and conceptual transformations underlying them. Scientometric data are therefore simultaneously a reduction of mathematical knowledge and a rich record of the scientific activity through which that knowledge is produced, organized, and transmitted. In this sense, the present analysis provides a coarse-grained view of the field that can serve as a starting point for more refined, concept-level investigations.

Seen from this perspective, bibliographic records are far from mere abstract counts~\cite{Aman2025}. They are traces left behind by a complex and evolving scientific ecosystem, capturing in fragmentary but systematic form some of the trajectories that have shaped the present organization of knowledge~\cite{Gibbons2010}. Such records constitute an important empirical basis for what has been called the \emph{computational history of science}~\cite{computational-history-2013}: much as archaeological records provide evidence for reconstructing past societies, bibliographic traces provide evidence for computationally investigating the development of scientific communities and their knowledge structures.

For category theory, this perspective offers both an external and an internal view of the field. By mapping its bibliographic traces, we provide a window for researchers interested in the evolution of mathematics and science and a mirror for the category theory community itself, making visible patterns, connections, and historical trends that might otherwise remain implicit or disappear from the community's collective view. In this way, scientometric analysis provides both a means of studying how mathematical communities emerge and evolve and a resource for collective self-understanding within the community.

In this work, we investigate how category theory interacts with other fields and how these interactions relate to the structure of its research community. Although a detailed account of knowledge flow would require semantic-level analyses of how results are incorporated into new research, keyword and citation patterns provide a useful large-scale proxy for identifying disciplinary connections and sources of influence.

\section{Our approach}
We conduct a large-scale data analysis of the zbMATH Open database\footnote{\url{https://zbmath.org/}} using publications explicitly classified under category theory (\emph{Mathematics Subject Classification} (MSC) code 18). Throughout the remainder of this article, whenever we refer to the \emph{category theory community}, we refer to this corpus of publications and its associated researchers, together with the various data structures and representations derived from it, such as co-authorship and topic networks. 

This corpus comprises publications within zbMATH that are classified under MSC 18. Consequently, it may underrepresent applications developed primarily outside mathematics, as well as mathematically related work that does not receive an MSC 18 classification. Nevertheless, it provides a useful basis for studying how category theory interacts with other mathematical areas and how those interactions have evolved over time.

First, we study the \emph{historical growth and thematic diversification} of the field, tracking changes in publication activity and in the distribution of Mathematical Subject Classification (MSC) topics over time. Second, we examine the \emph{structure of the research community and its patterns of scientific communication} through authorship, scientific output, collaboration, and the breadth of mathematical areas in which individual researchers work. Third, we investigate \emph{interactions across disciplinary boundaries} by analyzing the co-occurrence of MSC topics, the interfaces connecting category theory with other areas of mathematics and applied disciplines, and citation patterns between these domains as a proxy for knowledge flow.
These three perspectives provide a multi-layered view of the field, linking its conceptual development to observable changes in publication practices, collaboration networks, and the evolving language used to describe categorical research.

\subsection{Data Source}

Our primary data source is zbMATH Open, a comprehensive bibliographic database of mathematical publications curated by subject experts. Originally founded in 1931 as \emph{Zentralblatt für Mathematik und ihre Grenzgebiete}, the service has developed from a mathematical reviewing journal into a large-scale digital information system~\cite{Teschke2011}. Today, zbMATH Open is produced by the Berlin office of FIZ Karlsruhe and edited by the European Mathematical Society, FIZ Karlsruhe, and the Heidelberg Academy of Sciences and Humanities. The database currently indexes more than five million records across publication formats including journal articles, conference and collection articles, books, proceedings, dissertations, and preprints. Through the integration of the \emph{Jahrbuch über die Fortschritte der Mathematik}, it provides continuous coverage of the mathematical literature from 1868, together with earlier historical records. In 2021, zbMATH transitioned from a subscription-based service to the open-access platform zbMATH Open~\cite{Hulek2021,SchubotzTeschke2021}. The database provides structured metadata including publication year, authorship, publication venue, Mathematical Subject Classification (MSC) codes (version 2020), keywords, reviews and abstracts, references, and citation information~\cite{SchubotzTeschke2021}. Its long temporal coverage, systematic subject classification, and rich bibliographic metadata make zbMATH Open particularly suitable for bibliographical studies of the structure and development of mathematical research~\cite{HulekTeschke2024}. We use these features here to investigate the historical development, community structure, and disciplinary interactions of category theory.

\subsection{Data Extraction and Preprocessing}

We constructed a corpus of category theory publications by extracting records from zbMATH Open classified under MSC 18 (Category Theory), including both primary and secondary classifications. The resulting dataset comprises 33{,}947 documents published between 1942 and 2024. For each record, we retained available metadata including publication year, authorship, venue, MSC classifications, and citation information. Because some records have incomplete metadata, each analysis excludes those records that lack the relevant information. For instance, in the authorship analysis we excluded 776 documents that do not contain author information. As a result, absolute counts may differ slightly across analyses.
Basic preprocessing included the removal of duplicate entries and records with incomplete metadata. After these operations, the corpus used in the subsequent analyses comprised 30{,}854 records.

\subsection{Research Questions}
\label{subsec:research-questions}

Following the three aspects outlined above, we address the following research questions:

\begin{enumerate}
    \item How has the category theory research community developed over time in terms of publication activity and thematic diversification?
    
    \item What are the structural dynamics of the category theory research community, as reflected in patterns of authorship, scientific output, collaboration, and specialization?
    
    \item How does category theory interact with other areas of mathematics and related disciplines, and through which categorical interfaces are these interactions organized, as evidenced by topic co-occurrences and citation patterns?
    
\end{enumerate}

\section{Growth, Community Structure, and Knowledge Flow}

We organize our results around the three dimensions introduced above, corresponding to our three research questions: the historical growth and thematic diversification of category theory, the structure and communication patterns of its research community, and the disciplinary interactions through which categorical knowledge develops and circulates.

\subsection{Growth and Thematic Diversification}
\label{sec:growth-diversification}

\subsubsection{Growth and Expansion}

\begin{figure}[ht]
\centering
 \includegraphics[width=0.9\textwidth]{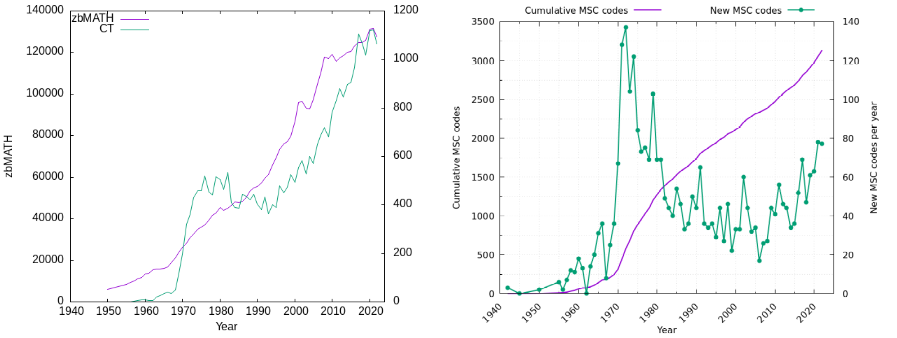}
\caption{Left: Annual number of publications in category theory and mathematics as a whole, showing their growth over time. Right: Number of newly observed MSC codes per year and cumulative number of distinct MSC codes associated with category theory publications (1942--2022). The yearly counts show the continued introduction of new topics, while the cumulative curve reflects the expansion of the field's topical scope across the Mathematics Subject Classification.}
    \label{fig:docs-year-msc-cum}
\end{figure}

Category theory has evolved into a highly active and steadily expanding field. After an initial phase of development, publication activity increased significantly in the 1970s, reaching approximately 500 papers per year. Following a period of slower growth in the 1980s, publication activity resumed its upward trend and accelerated after 2000 (Figure~\ref{fig:docs-year-msc-cum}).

In the past decade, annual publication counts have consistently exceeded 800. Moreover, over the last fifteen years, category theory has grown at a rate that outpaces mathematics as a whole. These trends indicate sustained publication growth, with a recent acceleration in activity within the field.

\subsubsection{Diversification of Topics}

The growth in publication activity is accompanied by an equally pronounced expansion in the range of topics associated with category theory. Using MSC codes as a proxy, we observe a continuous growth in the number of distinct topics associated with category theory publications\footnote{Some of the results in this section may be influenced by reclassifications of MSC codes (see subsection \ref{subsec:limitations}).}.

As shown by the cumulative curve in the right panel of Figure~\ref{fig:docs-year-msc-cum}, the number of distinct MSC codes associated with category theory grows from just a handful in the 1940s to over 3000 by 2022, with no indication of saturation. The introduction of new MSC codes, shown by the yearly-count curve in the same panel, accelerates through the early decades, peaking in the 1970s, and remains substantial thereafter. This suggests that category theory has steadily expanded its reach across the mathematical landscape.

\begin{figure}
    \centering
    \includegraphics[width=\linewidth]{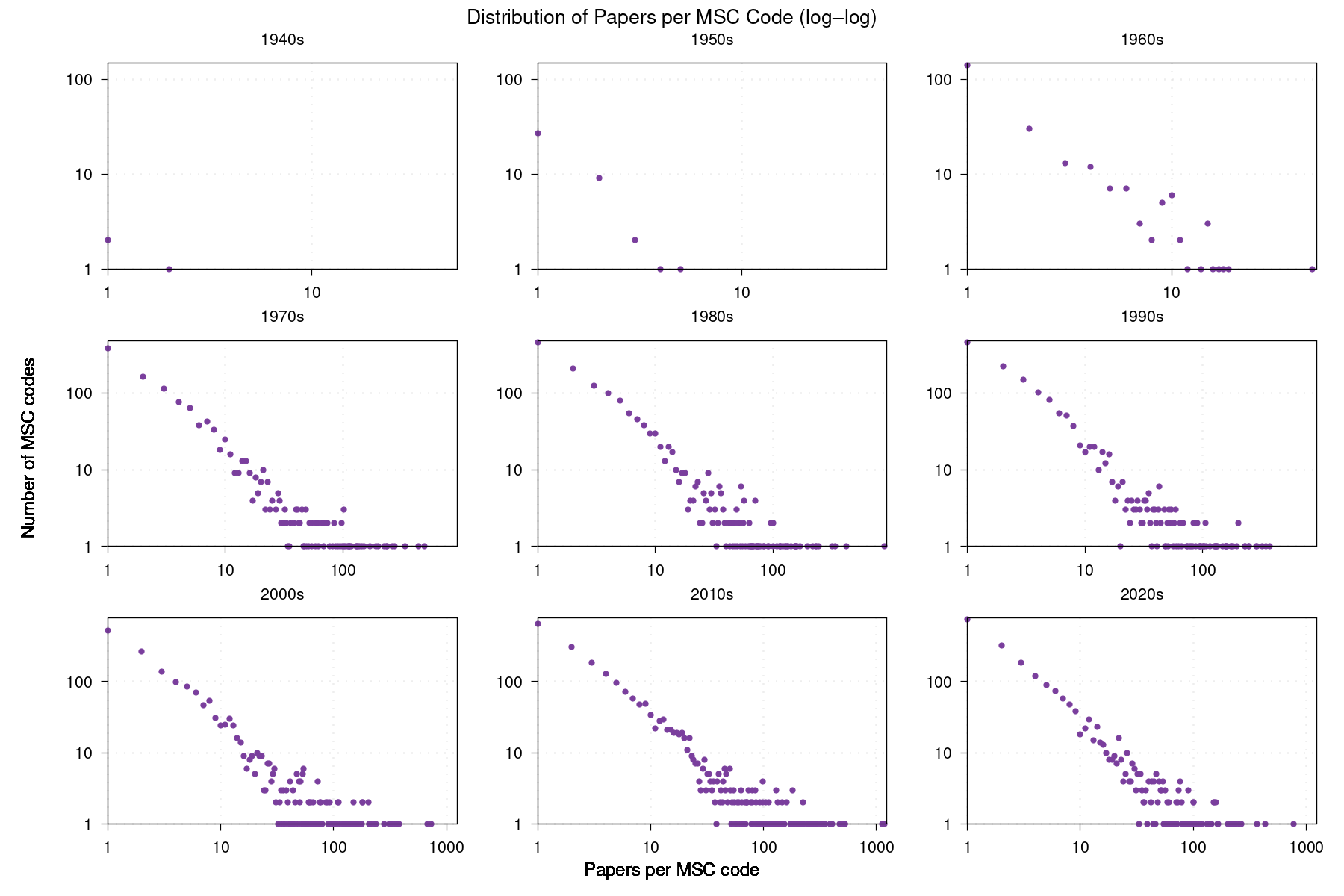}
    \caption{Log--log distributions of the number of papers per MSC code by decade (1940s--2020s), shown as nine panels (one per decade). In each panel, the horizontal axis represents the number of papers assigned to a given MSC code, and the vertical axis shows the number of MSC codes with that frequency. Each point corresponds to the count of codes appearing with a given frequency within the decade. 
    }
    \label{fig:topics-panels-decades}
\end{figure}

At the same time, the distribution of topics develops a pronounced heavy-tailed structure. As illustrated in Figure~\ref{fig:topics-panels-decades}, each decade exhibits a small number of highly represented MSC codes together with a long tail of infrequently occurring codes, a pattern that becomes increasingly pronounced over time. While a small number of MSC codes account for a significant fraction of publications, an increasing number of codes appear with low frequency.

\begin{figure}
    \centering
    \includegraphics[width=0.7\linewidth]{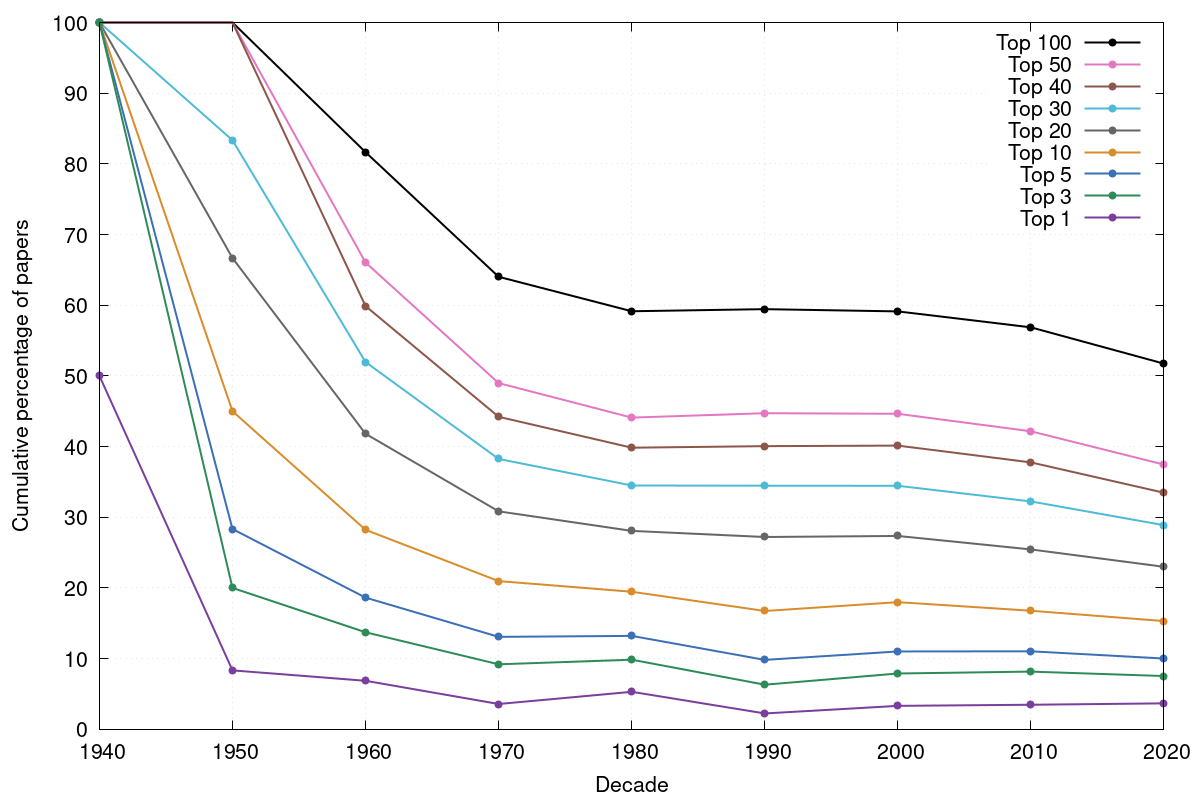}
    \caption{Cumulative share of publications accounted for by the top-$k$ most frequent MSC codes in each decade (1940s--2020s), for $k = 1, 3, 5, 10, 20, 30, 40, 50, 100$. Each curve shows the percentage of papers covered by the $k$ most frequently assigned MSC codes within a given decade. The results reveal a decreasing concentration over time: smaller sets of highly used codes account for a diminishing fraction of publications, indicating a progressive diversification of topics.}
    \label{fig:decade-topk-percent}
\end{figure}

Interestingly, over time, the dominance of the most frequent topics decreases. Figure~\ref{fig:decade-topk-percent} shows that the share of publications accounted for by the most common MSC codes steadily declines across decades: for example, the top 10 MSC codes account for roughly 45\% of publications in the 1950s, but only about 15--20\% in recent decades.

This combination of a persistent core and an expanding tail indicates a progressive diversification of the field: research activity is distributed across a growing range of topics, even as central themes remain identifiable. Together with the sustained increase in publication activity, these patterns show that the growth of category theory has been accompanied by a substantial expansion of its topical scope across the Mathematics Subject Classification.

\subsection{Community Structure and Scientific Communication}
\label{sec:community-structure}

\subsubsection{Community Growth and Scientific Output}

\begin{figure}
    \centering
    \includegraphics[width=\linewidth]{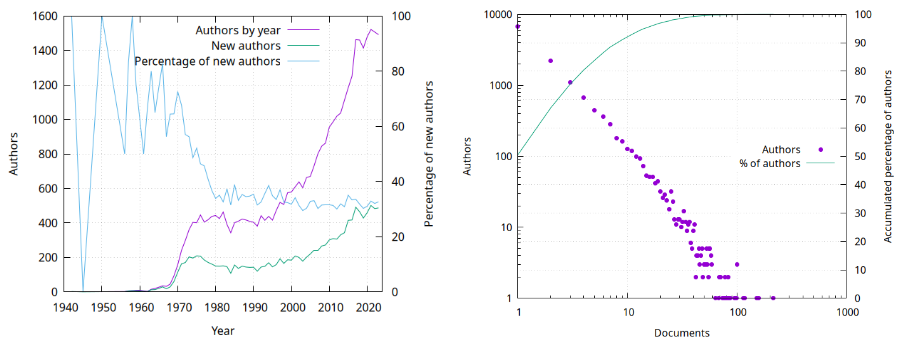}
    \caption{Left: Active authors by year. Number of authors that publish at least one document in a year (magenta). The percentage of new authors (blue) stabilized during the 80's. The number of new authors (green) has been increasing since circa 2000. Right: Distribution of the number of publications per author. Approximately 50\% of authors have contributed a single publication, whereas about 8\% have authored ten or more publications. The green curve shows the cumulative percentage of authors.}
    \label{fig:authors-year-docs}
\end{figure}

The growth of the field is mirrored by the expansion of its research community. Our dataset includes 13,392 authors, with the number of active authors per year increasing from a few dozen before 1960 to over 1,500 today (Figure~\ref{fig:authors-year-docs}, left).

While the absolute number of new authors continues to rise, their proportion has remained relatively stable since the 1980s (Figure~\ref{fig:authors-year-docs}, left). This indicates that the expansion of the community has combined a continuing influx of new researchers with an established core of researchers sustaining the field's development.

Authorship patterns also show substantial variation in publication frequency. Nearly half of all authors have one publication in the corpus, while around 8\% have 10 or more publications  (Figure~\ref{fig:authors-year-docs}, right). We found that only eight authors have contributed 100 documents or more. This distribution reflects a layered community, combining a broad base of occasional contributors with a stable core.

\subsubsection{Collaboration Patterns}
\begin{figure}
    \centering
    \includegraphics[width=.8\linewidth]{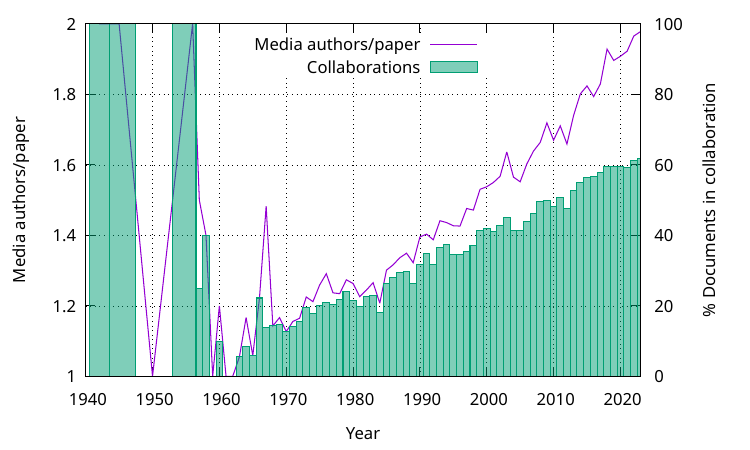}
    \caption{Evolution of scientific collaboration in the Category Theory corpus. The purple curve shows the annual median number of authors per publication, while the green boxes show the annual percentage of multi-author publications. Both measures exhibit a sustained upward trend over time. Right: Total number of publications as a function of the number of distinct collaborators. Authors with more collaborators generally have higher publication output.}
    \label{fig:collab-media}
\end{figure}

In scientometrics, collaboration is understood as document coauthorship.  Collaboration in science has increased in the last two decades \cite{Porter2009IsSB}, perhaps as a consequence of advancements in communication and technology \cite{Okamura2022AHO,Wagner2025InternationalSC}. Mathematicians are not the exception and their collaboration is increasing \cite{Grossman2002TheEO}, and lately has become even more notorious \cite{AgustnAquino2026SomePO}. The same pattern is observed in Category theory: since 1960, the media number of authors per paper has increased from 1.20 to 1.96 in 2024 (Figure \ref{fig:collab-media}).  Similarly, the proportion of multi-author publications has grown steadily, from 10\% in 1960 to 61.1\% in 2024.

\begin{figure}
    \centering
    \includegraphics[width=.8\linewidth]{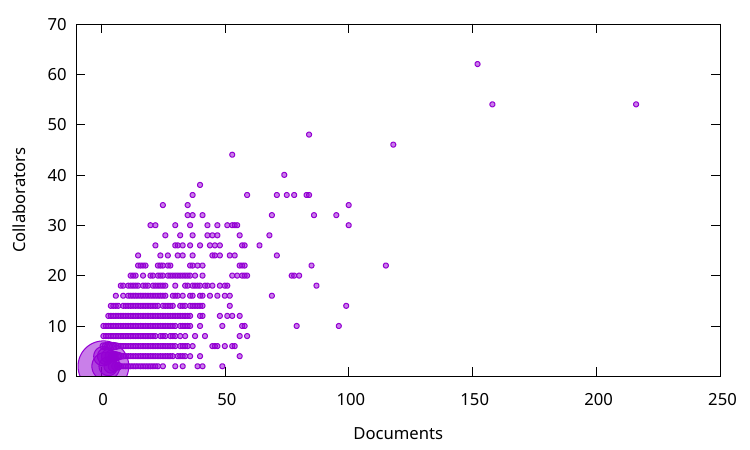}
    \caption{Total number of publications as a function of the number of distinct collaborators. Authors with more collaborators generally have higher publication output.}
    \label{fig:documents-vs-collaborators}
\end{figure}

Moreover, at the level of authors, collaboration is widespread: around 75\% of researchers have co-authored at least one paper, and authors with larger publication records also tend to have more co-authored work (Figure \ref{fig:documents-vs-collaborators}). Authors with 100 documents or more have from 29 to 61 collaborators (media of 41.1). The majority (56,6\%) of authors having 50 documents or more have more than 20 collaborators. Together with the increasing share of co-authored publications, these patterns show that collaboration has become an increasingly central mechanism for knowledge production in the community. 

\subsubsection{Propagation and External Reach}

\begin{figure}[ht]
    \centering
    \includegraphics[width=\linewidth]{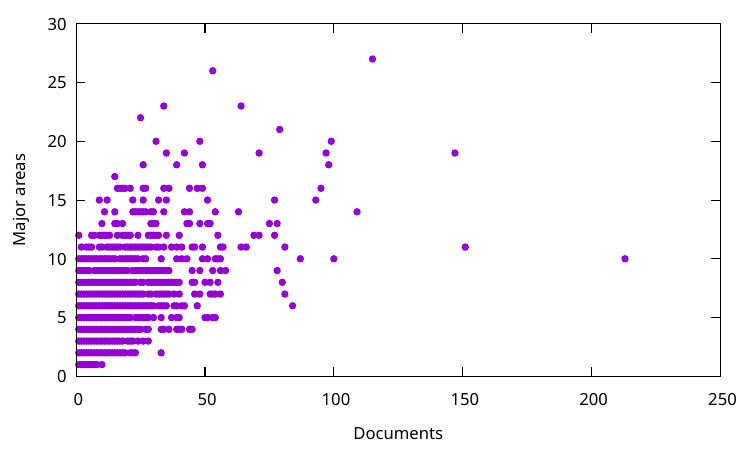}
    \caption{Number of major mathematical areas associated with each author as a function of publication count. A small subset of authors connects category theory with numerous MSC topics.}
    \label{fig:spec-major-area-ext}
\end{figure}

Beyond collaboration, researchers also differ substantially in the breadth of mathematical areas in which they work. While most authors in the category theory community publish within a relatively small number of major MSC areas, a small subset spans ten or more (Figure~\ref{fig:spec-major-area-ext}). A great majority of authors (96,4\%) have published documents on less than 5 areas, and only 5\% of prolific authors (published 50 documents or more) have descriptors on ten or more areas. Authors associated with ten or more areas are interdisciplinary researchers connecting category theory with multiple domains.

\begin{figure}
    \centering
        \includegraphics[width=\linewidth]{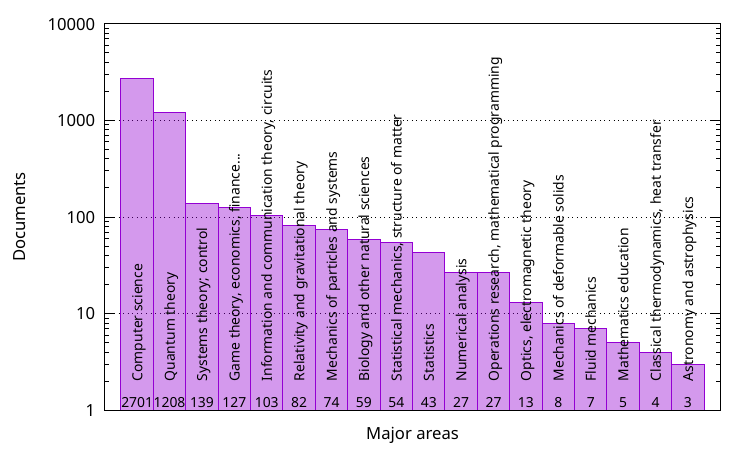}
    \caption{Frequency of major areas in the extension set (documents having at least one MSC code within 60-97), showing that computer science (MSC: 68) and physics (MSC: 70-86) are the most frequent application domains of category theory.}
    \label{fig:mayor-areas-vs-docs}
\end{figure}

 As shown in Figure~\ref{fig:mayor-areas-vs-docs}, computer science (MSC: 68) is by far the most frequent external domain with 2701 documents, followed by physics (MSC: 70-86) with 1378. Beyond these two dominant areas, the distribution extends across systems theory (MSC: 93) with 139, Game theory and economics (MSC: 91) with 127, Information and communication (MSC: 94) with 103,  biology (MSC: 92) with 59, and several other applied domains, covering a broader and more heterogeneous landscape of applications.

These author-level patterns show that some researchers contribute to category theory publications spanning a particularly broad range of mathematical and applied areas, potentially acting as \emph{propagators} of categorical ideas across disciplinary boundaries. The next subsection shifts from researchers to topics, examining the disciplinary interactions through which these connections are organized and identifying the specific interfaces linking category theory with other areas of mathematics and its applications.

\subsection{Disciplinary Interactions and Knowledge Flow}

We now move from the temporal and author-level structure of the field to the disciplinary interactions through which categorical research connects with other domains. We examine these interactions at three complementary levels. First, we analyze co-occurrences among MSC topics to identify the principal interfaces between category theory, other areas of mathematics, and applied disciplines (Section~\ref{sec:topic-interaction}). We then classify publications according to whether their topical composition remains within category theory, extends to neighboring mathematical areas, or reaches applied domains (Section~\ref{sec:core-adjacency-extension}). Finally, we use citation patterns to examine the flow of information between these regions and, in particular, the permeability of the category-theoretic core (Section~\ref{sec:core-is-permeable}).

\subsubsection{Interactions Between Topics}\label{sec:topic-interaction} 

We investigate how category theory interacts internally and with other mathematical and applied domains. To this end, the list of topics was constructed using three-digit MSC codes (for example, 18A—General theory of categories and functors). Based on this classification, we identify as \emph{Category Theory} all topics whose MSC codes begin with 18. The remaining topics are grouped into two classes: \emph{Applied Topics}, corresponding to MSC codes 60–97, and \emph{Other Topics}, which comprises all remaining MSC codes.

We represent these interactions using \emph{chord diagrams}, in which each chord represents the frequency of co-occurrence between two three-digit MSC descriptors. Each article contributes one co-occurrence to every pair of descriptors it contains, so the thickness of a chord reflects the strength of the interaction between the corresponding topics. Because articles are often assigned multiple MSC descriptors, each publication generally contributes multiple topic pairs, so co-occurrence counts should not be interpreted as numbers of publications.  Moreover, since MSC classifications are aggregated to the three-digit level, distinct four- or five-digit descriptors may collapse to the same three-digit descriptor. Consequently, an article may contribute co-occurrences between descendants of the same three-digit MSC descriptor, which appear as self-loops in the chord diagrams.

\begin{figure}
    \centering
    \includegraphics[width=\linewidth]{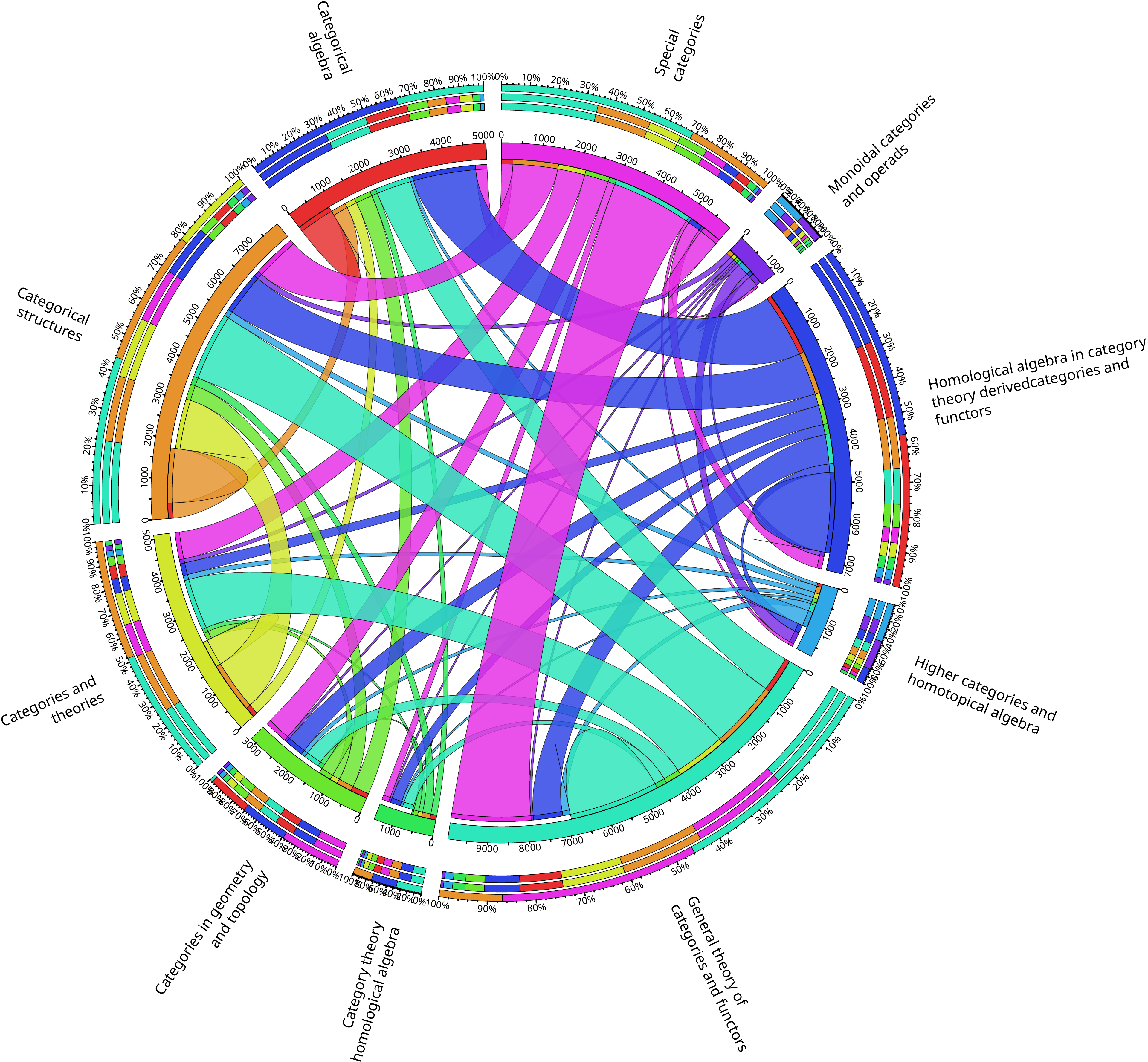}
    \caption{Chord diagram of co-occurences between Category Theory topics in the zbMath Open corpus of 33,947 Category Theory publications (MSC 18, primary or secondary classification, 1942--2024).  Nodes correspond to three-digit MSC descriptors beginning with 18, and each chord represent the number of articles simultaneously classified under the corresponding pair of descriptors.  The thickness of each chord is proportional to the frequency of co-occurence. Because MSC classifications are aggregated to the three-digit level, distinct more specific descriptors belonging to the same three-digit topic may produce self-loops. }
    \label{fig:ct-ct}
\end{figure}

\begin{figure}
    \centering
    \includegraphics[width=\linewidth]{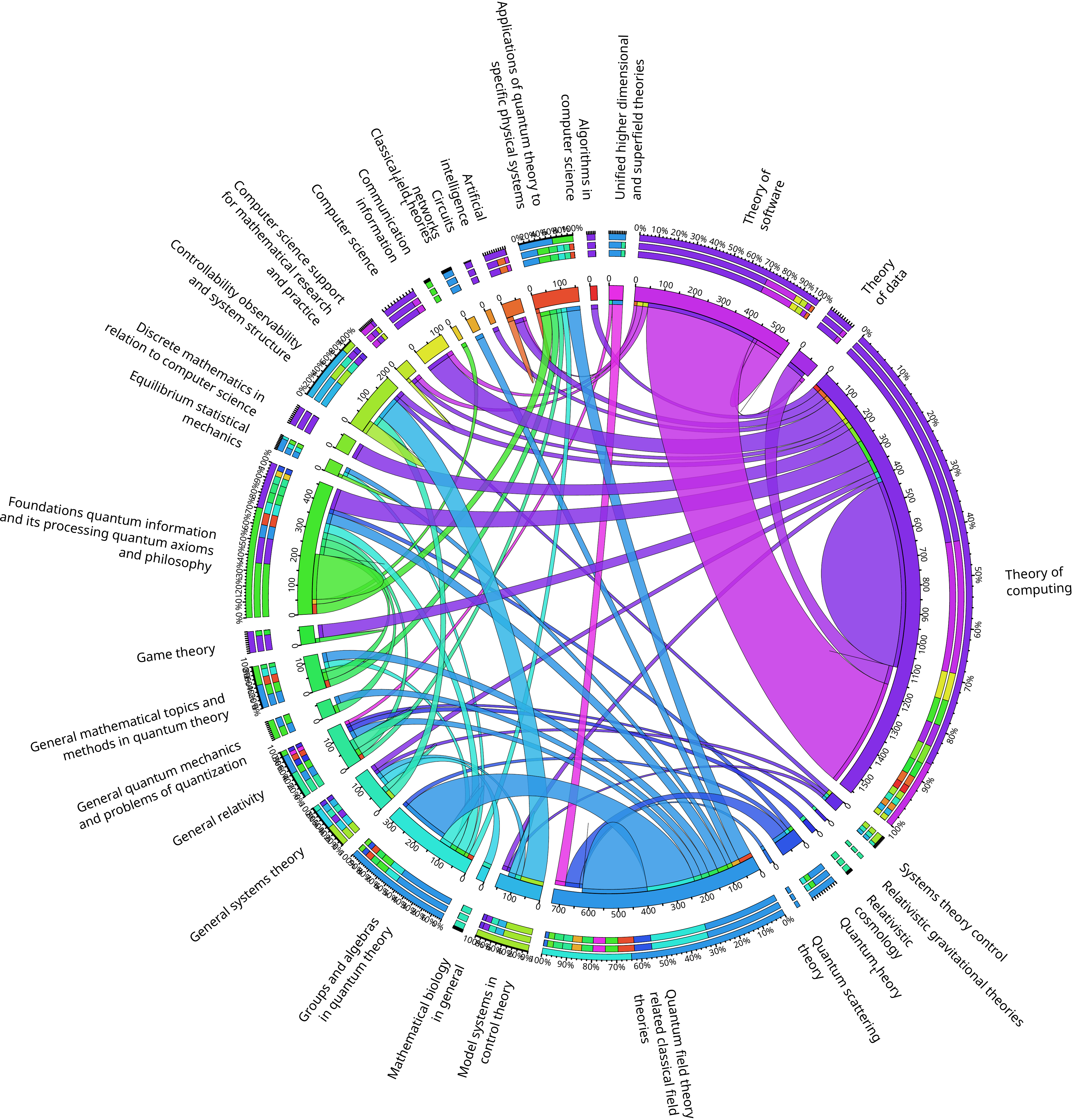}
    \caption{Chord diagram of co-occurrences between Applied Topics in the zbMATH Open corpus of 33,947 Category Theory publications (MSC 18, primary or secondary classification, 1942--2024). Nodes correspond to three-digit MSC descriptors in the Applied Topics class (MSC 60--97), and each chord represents the number of articles simultaneously classified under the corresponding pair of descriptors. The thickness of each chord is proportional to the frequency of co-occurrence. Because MSC classifications are aggregated to the three-digit level, distinct more specific descriptors belonging to the same three-digit topic may produce self-loops. }
    \label{fig:app-app}
\end{figure}

\begin{figure}
    \centering
    \includegraphics[width=\linewidth]{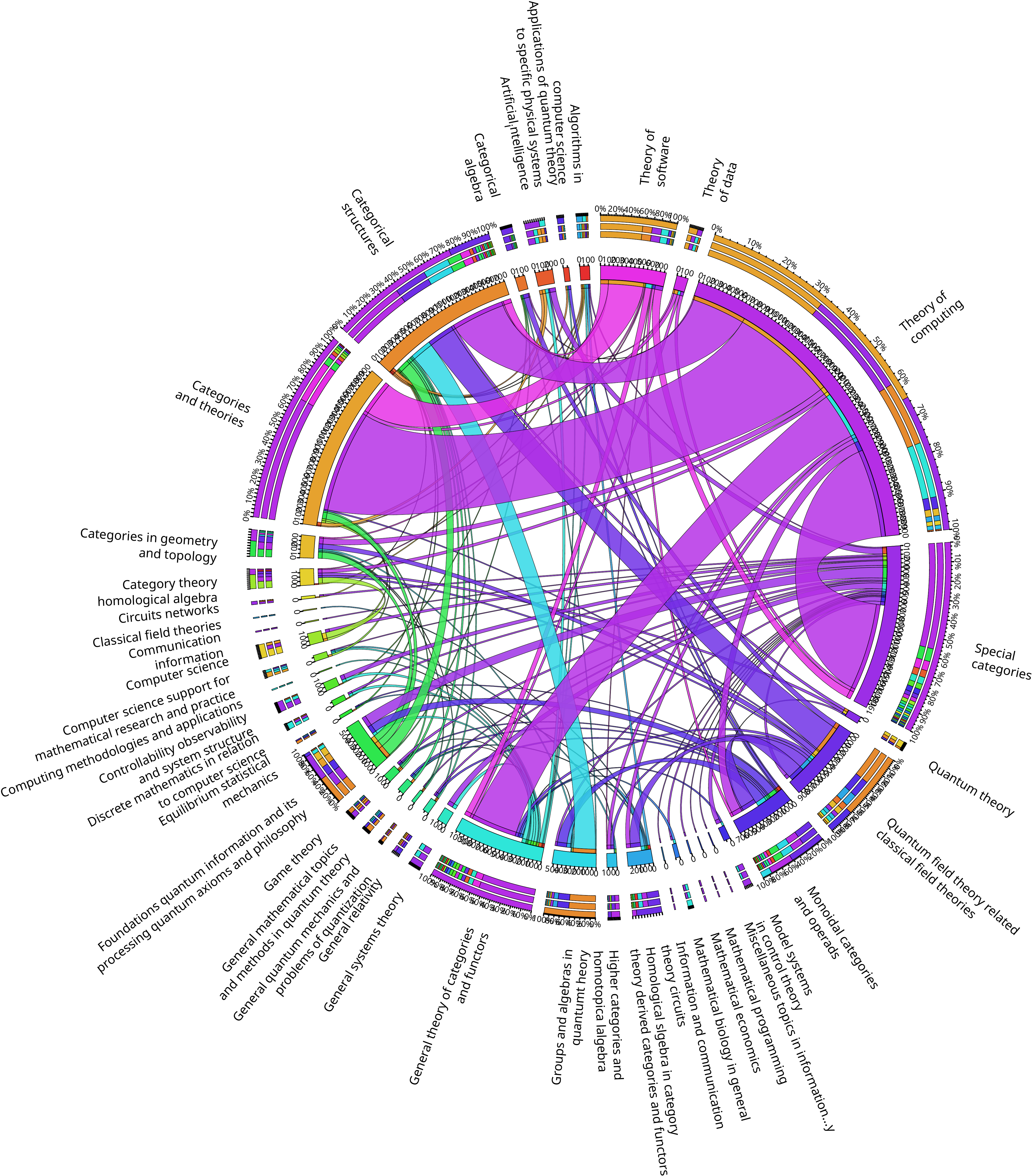}
    \caption{Chord diagram of co-occurrences between Category Theory topics and Applied Topics in the zbMATH Open corpus of 33,947 Category Theory publications (MSC 18, primary or secondary classification, 1942--2024). Nodes correspond to three-digit MSC descriptors in Category Theory (MSC 18) and Applied Topics (MSC 60--97), and each chord represents the number of articles simultaneously classified under the corresponding pair of descriptors, with one descriptor from each class. The thickness of each chord is proportional to the frequency of co-occurrence.}
    \label{fig:app-ct}
\end{figure}

\begin{figure}
    \centering
    \includegraphics[width=\linewidth]{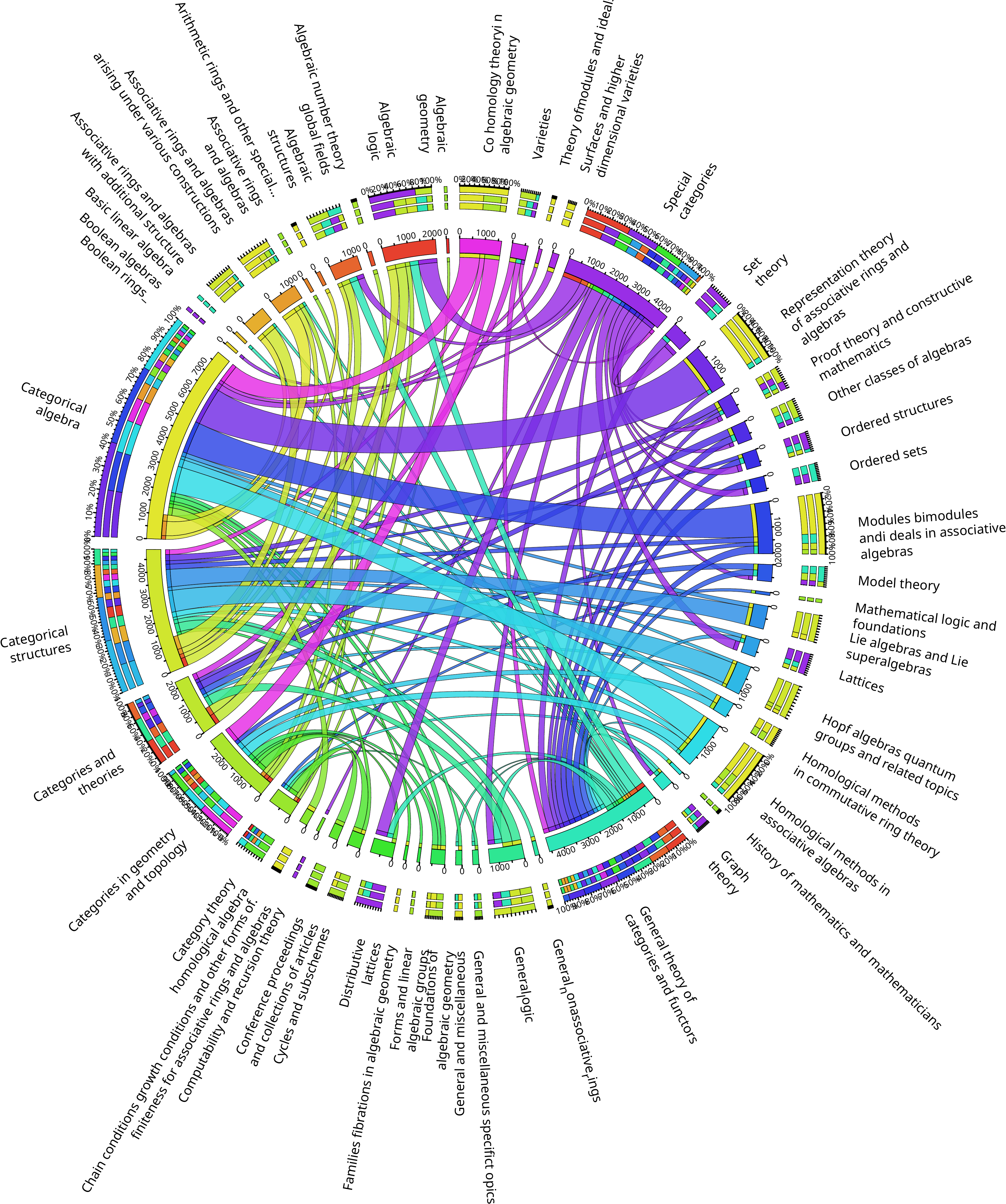}
    \caption{Chord diagram of co-occurrences between Category Theory topics and Other Topics in the zbMATH Open corpus of 33,947 Category Theory publications (MSC 18, primary or secondary classification, 1942--2024). Nodes correspond to three-digit MSC descriptors in Category Theory (MSC 18) and Other Topics (MSC codes below 60, excluding MSC 18), and each chord represents the number of articles simultaneously classified under the corresponding pair of descriptors, with one descriptor from each class. The thickness of each chord is proportional to the frequency of co-occurrence.}
    \label{fig:ct-others}
\end{figure}

The three classes of MSC topics yield six possible types of interaction: three within classes and three between classes. We represent each interaction type with a separate chord diagram. Interactions within Category Theory are shown in Figure \ref{fig:ct-ct}, within Applied Topics in Figure~\ref{fig:app-app}, and within Other Topics in Figure~\ref{fig:others-others}.  Interactions between Applied Topics and Category Theory are shown in Figure~\ref{fig:app-ct}, between Category Theory and Other Topics in Figure~\ref{fig:ct-others}, and between Applied Topics and Other Topics in  Figure~\ref{fig:app-others}.

We first examine the internal interaction structure of each of the three classes before turning to the interfaces between them.

\paragraph{Interactions within Category Theory topics}
Figure~\ref{fig:ct-ct} shows two distinct interaction patterns among the MSC descriptors of Category Theory. A first group, consisting of \emph{General theory of categories and functors}, \emph{Special categories}, \emph{Categories and theories}, \emph{Categorical structures}, \emph{Categorical algebra}, and \emph{Homological algebra in category theory, derived categories and functors}, forms the most strongly interconnected part of the diagram. In contrast, \emph{Categories in geometry and topology}, \emph{Category theory homological algebra}, \emph{Monoidal categories and operads}, and \emph{Higher categories and homotopical algebra} exhibit substantially weaker interactions with the rest of the network.

Within Category Theory, the strongest interactions occur between \emph{General theory of categories and functors} and \emph{Special categories} (2038 co-occurrences), \emph{General theory of categories and functors} and \emph{Categorical structures} (1733), \emph{Categorical algebra} and \emph{Homological algebra in category theory, derived categories and functors} (1626), \emph{General theory of categories and functors} and \emph{Categories and theories} (1318), and \emph{Categorical structures} and \emph{Categories and theories} (1280).

Among the more specialized descriptors, \emph{Categories in geometry and topology} connects most strongly with \emph{Special categories} (617), \emph{Homological algebra in category theory} (491), \emph{Categorical algebra} (434), and \emph{General theory of categories and functors} (409). Meanwhile, \emph{Monoidal categories and operads} exhibits its strongest interaction with \emph{Higher categories and homotopical algebra} (306), followed by \emph{Categorical structures} (142).

\begin{figure}
    \centering
    \includegraphics[width=\linewidth]{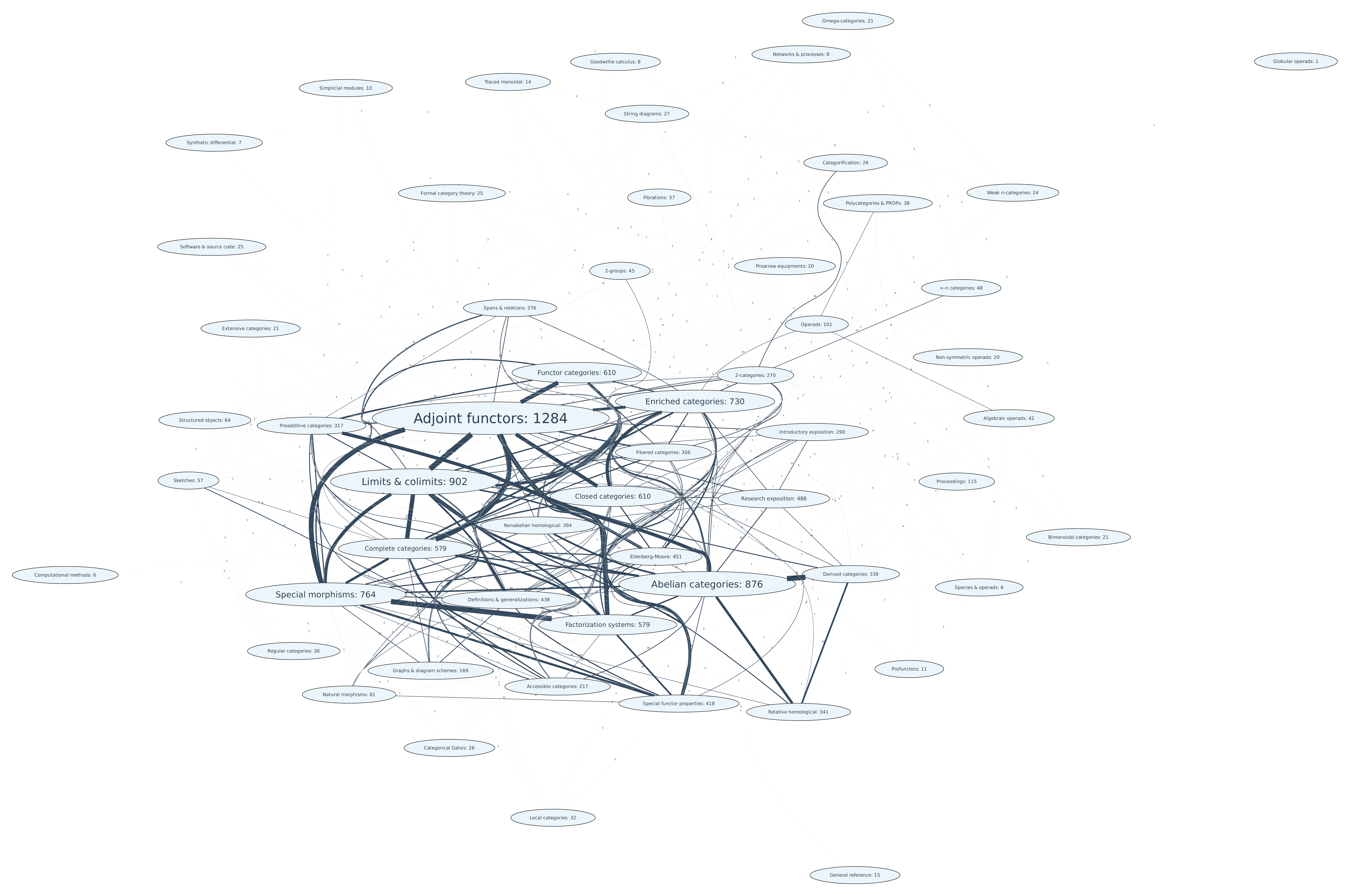}
    \caption{Network of Category Theory topics at the five-digit MSC level. Undirected network of co-occurrences between Category Theory topics in the zbMATH Open corpus of 33,947 Category Theory publications (MSC 18, primary or secondary classification, 1942--2024). Nodes correspond to five-digit MSC descriptors; topic names have been shortened for readability. Each edge connects two descriptors assigned to the same publication; its weight is the number of publications in which the pair co-occurs. Edge thickness is proportional to this co-occurrence frequency, while node size is proportional to weighted degree, defined here as the total frequency of co-occurrences of the corresponding topic with all other Category Theory topics. There are few MSC descriptors concerning the type of publication rather than topics (e.g., proceedings, introductory expositions, among others), which we included for the sake of completeness in our analysis.}
    \label{fig:ct-ct-zoomed-in}
\end{figure}

The three-digit analysis shows the large-scale organization of Category Theory topics, but its coarse resolution aggregates many mathematically distinct constructions within the same MSC topic. We therefore repeat the analysis at the five-digit level to examine the internal structure underlying these aggregate interactions. Figure~\ref{fig:ct-ct-zoomed-in} shows the resulting network.

At this finer resolution, the strong connectivity observed for \emph{General theory of categories and functors} (18A) resolves into interactions among specific categorical constructions. Several of the strongest pairwise interactions in the entire network occur within this class (Table~\ref{tab:ct-five-digit-interactions}), including \emph{Limits and colimits}--\emph{Adjoint functors} (147 co-occurrences), \emph{Special morphisms}--\emph{Factorization systems} (138), \emph{Special morphisms}--\emph{Adjoint functors} (131), \emph{Complete categories}--\emph{Adjoint functors} (126), and \emph{Factorization systems}--\emph{Adjoint functors} (120). Among these, \emph{Adjoint functors} (18A40) emerge as a particularly broadly connected topic. In addition to several of the strongest interactions within 18A, adjunctions connect to topics belonging to other three-digit classes, including \emph{Eilenberg--Moore constructions} (18C20), \emph{Closed categories} (18D15), \emph{Enriched categories} (18D20), and \emph{Abelian categories} (18E10).

The finer network also clarifies other aggregate patterns observed in Figure~\ref{fig:ct-ct}. In particular, the strong interaction between \emph{Categorical algebra} (18E) and \emph{Homological algebra in category theory, derived categories and functors} (18G) is resolved into several prominent relationships. The interaction between \emph{Abelian categories} (18E10) and \emph{Derived categories} (18G80), with 140 co-occurrences, is the second strongest pairwise interaction in the entire five-digit network. \emph{Abelian categories} also interacts strongly with \emph{Preadditive categories} (18E05; 91), \emph{Relative homological algebra} (18G25; 67), and \emph{Nonabelian homological algebra} (18G50; 48). Likewise, the connectivity of \emph{Categorical structures} (18D) can be traced in part to \emph{Closed categories} and \emph{Enriched categories}, which co-occur 89 times and also connect strongly with the general categorical constructions centered around 18A.

\begin{table}[ht]
\centering
\caption{Selected strong co-occurrences between Category Theory topics at the five-digit MSC level, ordered by decreasing co-occurrence frequency. Each pair represents two MSC topics assigned to the same publication, and the co-occurrence count gives the number of publications classified under both topics.}
\label{tab:ct-five-digit-interactions}
\begin{tabular}{llr}
\toprule
\textbf{Topic 1} & \textbf{Topic 2} & \textbf{Co-occurrences} \\
\midrule
18A30 --- Limits and colimits
& 18A40 --- Adjoint functors
& 147 \\
18E10 --- Abelian categories
& 18G80 --- Derived categories
& 140 \\
18A20 --- Special morphisms
& 18A32 --- Factorization systems
& 138 \\
18A20 --- Special morphisms
& 18A40 --- Adjoint functors
& 131 \\
18A35 --- Complete categories
& 18A40 --- Adjoint functors
& 126 \\
18A32 --- Factorization systems
& 18A40 --- Adjoint functors
& 120 \\
18A30 --- Limits and colimits
& 18A35 --- Complete categories
& 117 \\
18A25 --- Functor categories
& 18A40 --- Adjoint functors
& 113 \\
18A40 --- Adjoint functors
& 18D15 --- Closed categories
& 99 \\
18E05 --- Preadditive categories
& 18E10 --- Abelian categories
& 91 \\
18A20 --- Special morphisms
& 18A30 --- Limits and colimits
& 90 \\
18D15 --- Closed categories
& 18D20 --- Enriched categories
& 89 \\
\bottomrule
\end{tabular}
\end{table}

\paragraph{Interactions within Applied Topics}

The interaction network between Applied Topics exhibits three main centers. The largest is a computer science cluster centered on \emph{Theory of computing} and \emph{Theory of software}. Besides the self-interaction of \emph{Theory of computing} (673), the strongest interaction in the network is between these two descriptors (397). The cluster further includes strong interactions with \emph{Computer Science} (84), \emph{Foundations of quantum information} (77), \emph{Theory of data} (67), \emph{Discrete mathematics in relation to computer science} (57), and \emph{Game Theory} (47), making theoretical computer science the principal interface through which category theory connects to a broad range of applied disciplines.

A second center is formed by quantum theory, where \emph{Quantum field theory and related classical field theories} interacts most strongly with \emph{Groups and algebras in quantum theory} (171). Additional strong interactions include \emph{Quantum theory} (52), \emph{Applications of quantum theory to specific physical systems} (49), \emph{Foundations of quantum information} (38), and \emph{Unified, higher-dimensional, and super field theories} (37). Notably, \emph{Foundations of quantum information} provides the principal bridge between the computer science and quantum clusters through its interaction with \emph{Theory of computing}.

A third center is formed by systems and control, organized around \emph{Controllability, observability, and system structure}, \emph{Model systems in control theory}, and \emph{General systems theory}. Its strongest interaction is between the first two descriptors (83), followed by the interaction between \emph{General systems theory} and \emph{Controllability, observability, and system structure} (32). Unlike the previous two clusters, these interactions are largely internal, with comparatively few interactions linking systems and control to other application areas. Biological applications likewise remain peripheral in the network. Their only notable interaction is between \emph{General systems theory} and \emph{Mathematical biology in general} (29).

\emph{Artificial intelligence} (68T) remains a comparatively specialized application area. Beyond 26 self-co-occurrences, it interacts primarily with \emph{Theory of computing} (68Q) (32) and, to a lesser extent, with \emph{Theory of software} (68N) (11), exhibiting essentially no direct interactions with the remaining applied topics. 

Having characterized the internal organization of the three topic classes (see Supplementary Material \ref{supplement}), we now turn to interactions between them.

\paragraph{Interactions of Category Theory and Applied Topics}
One interesting finding is that category theory interacts with application domains through different interfaces. One possible interpretation is that different application domains encode different forms of compositionality and therefore engage with different category-theoretic ideas and constructions.

The interactions between Category Theory and Applied Topics are dominated by theoretical computer science, followed by quantum theory, with systems and control forming a smaller but clearly identifiable interface. 

The strongest interactions between category theory and computer science all involve \emph{Theory of computing}. These are, in decreasing order, its interactions with \emph{Categories and theories} (1285 co-occurrences), \emph{Special categories} (1092), \emph{Categorical structures} (664), and \emph{General theory of categories and functors} (573). Together with the interaction between \emph{Categories and theories} and \emph{Theory of software} (434), these results identify theoretical computer science as the principal applied interface of category theory.

The second major interface is quantum theory. Here, \emph{Quantum field theory and related classical field theories} interacts most strongly with \emph{Categorical structures} (337), followed by \emph{Groups and algebras in quantum theory} (273), \emph{Monoidal categories and operads} (183), and \emph{Homological algebra in category theory} (127). 

A notable feature of the diagram is that theoretical computer science and quantum theory connect to category theory through distinct categorical interfaces. Whereas the computer science interface is concentrated around the foundational topics \emph{Categories and theories}, \emph{Special categories}, and \emph{General theory of categories and functors}, the quantum interface is centered on \emph{Categorical structures}, \emph{Monoidal categories and operads}, and \emph{Homological algebra in category theory}. 

Systems and control exhibit a smaller but coherent interaction pattern. The strongest connections are between \emph{Special categories} and \emph{General systems theory} (72), \emph{Controllability, observability, and system structure} (57), and \emph{Model systems in control theory} (56), indicating a comparatively narrow interface concentrated around a single category-theoretic topic.

Beyond these dominant interfaces, several additional application areas emerge. Artificial intelligence interacts primarily with \emph{Special categories} (74 co-occurrences), followed by \emph{General theory of categories and functors} (41), \emph{Categories and theories} (37), and \emph{Categorical structures} (36). Table~\ref{tab:msc68} shows that, among the MSC~68 three-digit descriptors, interactions with category theory are overwhelmingly concentrated in theoretical computer science. In particular, \emph{Theory of computing}, \emph{Theory of software}, and \emph{Theory of data} account for 3962, 774, and 136 total co-occurrences, respectively, whereas the other MSC~68 descriptor contributes 204, 115, 73, and 12 interactions. 

Other application areas remain comparatively small. \emph{Circuits and networks} is connected principally to \emph{Special categories} (21) and \emph{Categorical structures} (12), revealing a distinct but limited interface with network theory. Likewise, \emph{Game theory, economics, social and behavioral sciences} interacts mainly with \emph{Special categories} (28), \emph{Categorical structures} (24), \emph{Categories and theories} (19), \emph{Monoidal categories and operads} (16), and \emph{General theory of categories and functors} (13). Finally, \emph{Mathematical programming} exhibits a single notable interaction with \emph{Special categories} (11).

\begin{table}[ht]
\centering
\caption{Total co-occurrences between each MSC~68 descriptor and Category Theory topics.}
\label{tab:msc68}
\begin{tabular}{lr}
\toprule
\textbf{MSC 68 descriptor} & \textbf{Total interactions} \\
\midrule
Theory of computing (68Q) & 3962 \\
Theory of software (68N) & 774 \\
Artificial intelligence (68T) & 204 \\
Theory of data (68P) & 136 \\
Algorithms in computer science (68W) & 115 \\
Computer science support for mathematical research and practice (68V) & 73 \\
Computing methodologies and applications (68U) & 12 \\
\bottomrule
\end{tabular}
\end{table}

\paragraph{Interactions of Category Theory and Other Topics}

Figure~\ref{fig:ct-others} shows that different branches of mathematics (MSC codes within Other Topics) connect to category theory through distinct category theory topics rather than a common categorical interface.

The principal interface is \emph{Homological algebra in category theory, derived categories and functors}. It accounts for the strongest interaction in the entire diagram with \emph{Homological methods in associative algebras} (2269 co-occurrences), followed by \emph{Applied homological algebra and category theory in algebraic topology} (1800), \emph{Homological methods in commutative ring theory} (1095), and \emph{Homotopy theory} (1081).

A second major interface is \emph{Categorical algebra}, whose strongest interactions are with \emph{Representation theory of associative rings and algebras} (1598), \emph{Modules, bimodules and ideals in associative algebras} (1351), \emph{Homological methods in associative algebras} (1111), and \emph{Cohomology theory in algebraic geometry} (821).

Other branches of category theory show similarly distinctive interaction profiles. \emph{Special categories} interacts primarily with general topology, including \emph{Basic constructions in general topology} (1164), \emph{Generalities in topology} (632), and \emph{Fairly general properties of topological spaces} (508), while also exhibiting strong interactions with \emph{Algebraic logic} (778) and \emph{Set theory} (536). In contrast, \emph{Categorical structures} is linked most strongly to \emph{Quantum groups and related algebraic methods} (897), \emph{Homotopy theory} (823), \emph{Categories in geometry and topology} (721), and \emph{Lie algebras and Lie superalgebras} (734).

For interaction of Applied Topics and Other Topics see Supplementary Material \ref{supplement}.

\subsubsection{Core, Adjacency, and Extension}\label{sec:core-adjacency-extension} 
In Section~\ref{sec:topic-interaction}, we classified three-digit MSC topics into three groups---Category Theory, Other Topics, and Applied Topics---and analyzed their co-occurrence. We now shift our analysis from topics to publications, classifying each document according to the three topic groups associated with its MSC descriptors. We distinguish three classes: \emph{core}, comprising publications whose descriptors all belong to Category Theory (MSC 18); \emph{adjacency}, comprising publications with descriptors from Category Theory and Other Topics, but none from Applied Topics; and \emph{extension}, comprising publications with at least one descriptor from Applied Topics.

The majority of publications belong to the adjacency set (72.7\%), while 14.2\% belong to the core and 13.1\% to the extension set (Table~\ref{tab:docs-per-set}). This indicates that category theory is primarily developed in interaction with other areas of mathematics, rather than in isolation.

\begin{table}[h]
    \centering
    \caption{Documents per set}
    \label{tab:docs-per-set}
    \begin{tabular}{ccc}
        \hline
        Core & Adjacency & Extension \\
        \hline
        4664 & 23{,}806 & 4{,}287 \\
        14.2\% & 72.7\% & 13.1\% \\
        \hline
    \end{tabular}
\end{table}

\begin{figure}
    \centering
    \includegraphics[width=\linewidth]{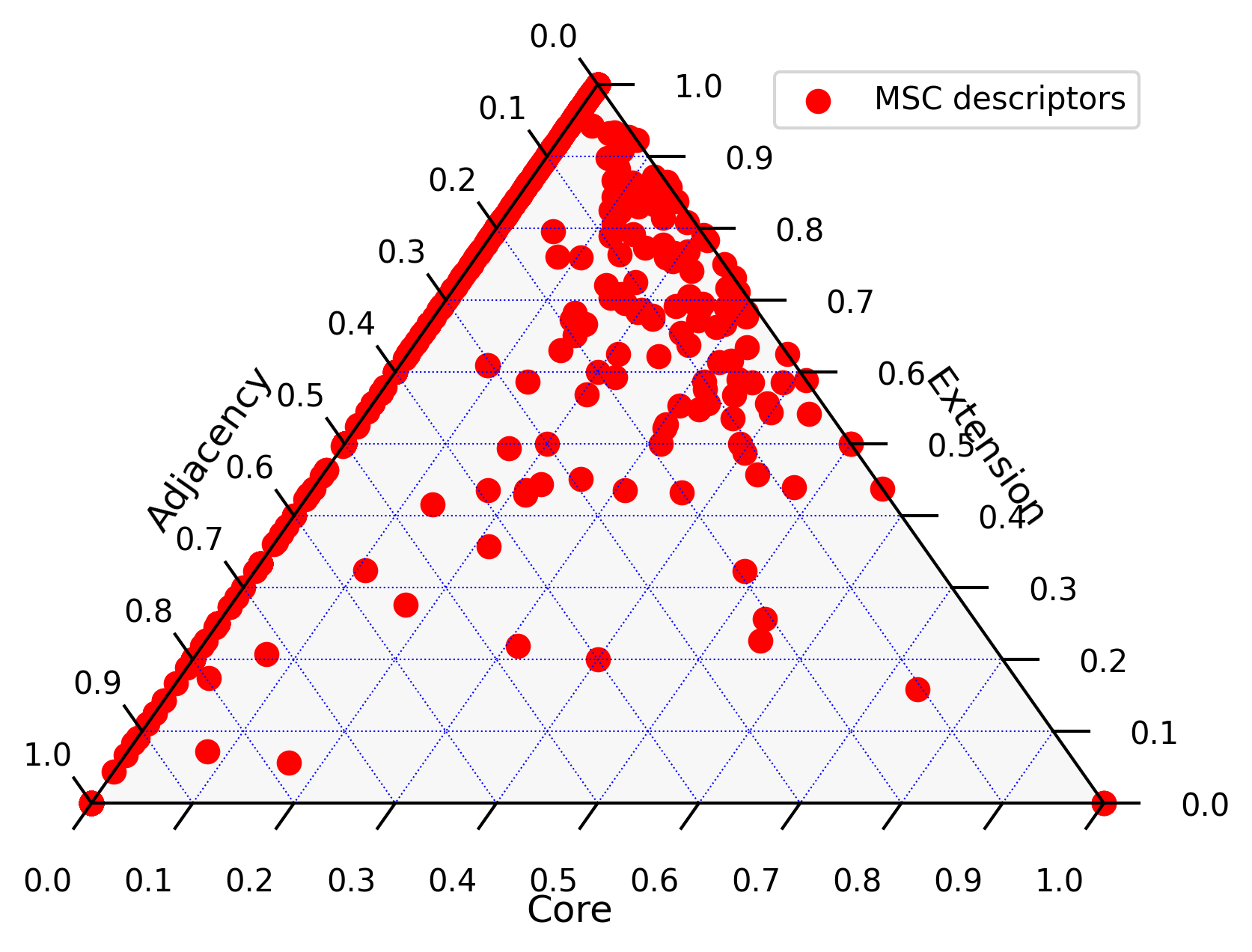}
    \caption{Relative frequency of descriptors. Each point represents a MSC descriptor and its position represent the relative frequency of documents containing it that fall in each set. }
    \label{fig:descriptors2}
\end{figure}

Category Theory descriptors are not confined to the core set, but also occur throughout the adjacency and extension sets (Figure~\ref{fig:descriptors2}). More importantly, the adjacency set is concentrated around algebra, topology, and logic, recovering at the publication level the principal neighboring mathematical areas identified in the topic-level interaction analysis of Section~\ref{sec:topic-interaction}.

\begin{figure}
    \centering
    \includegraphics[width=\linewidth]{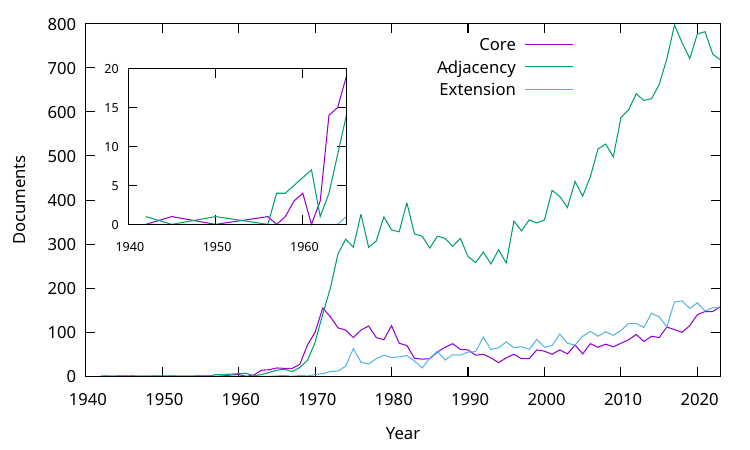}
    \caption{Number of documents per year. Core corresponds to documents having only MSC codes of Category Theory; extension corresponds to document having at least one MSC between 60-97; adjacency to anything else.}
    \label{fig:descriptors}
\end{figure}

The relative prevalence of the three sets has also changed over time (Figure~\ref{fig:descriptors}). Core and adjacency publications occurred at comparable levels during the early decades, with neither consistently dominating the other. From the 1970s onward, adjacency became increasingly dominant, although both core and adjacency declined until around 2000. Over much of this period, extension publications followed a different trajectory, increasing while the other two sets declined and surpassing core publications around the 1990s. Since the 2000s, all three sets have resumed growth, with adjacency remaining the largest.

\subsubsection{A Permeable Category-Theoretic Core}\label{sec:core-is-permeable} 
\begin{figure}[h!]
    \centering
    \includegraphics[width=\linewidth]{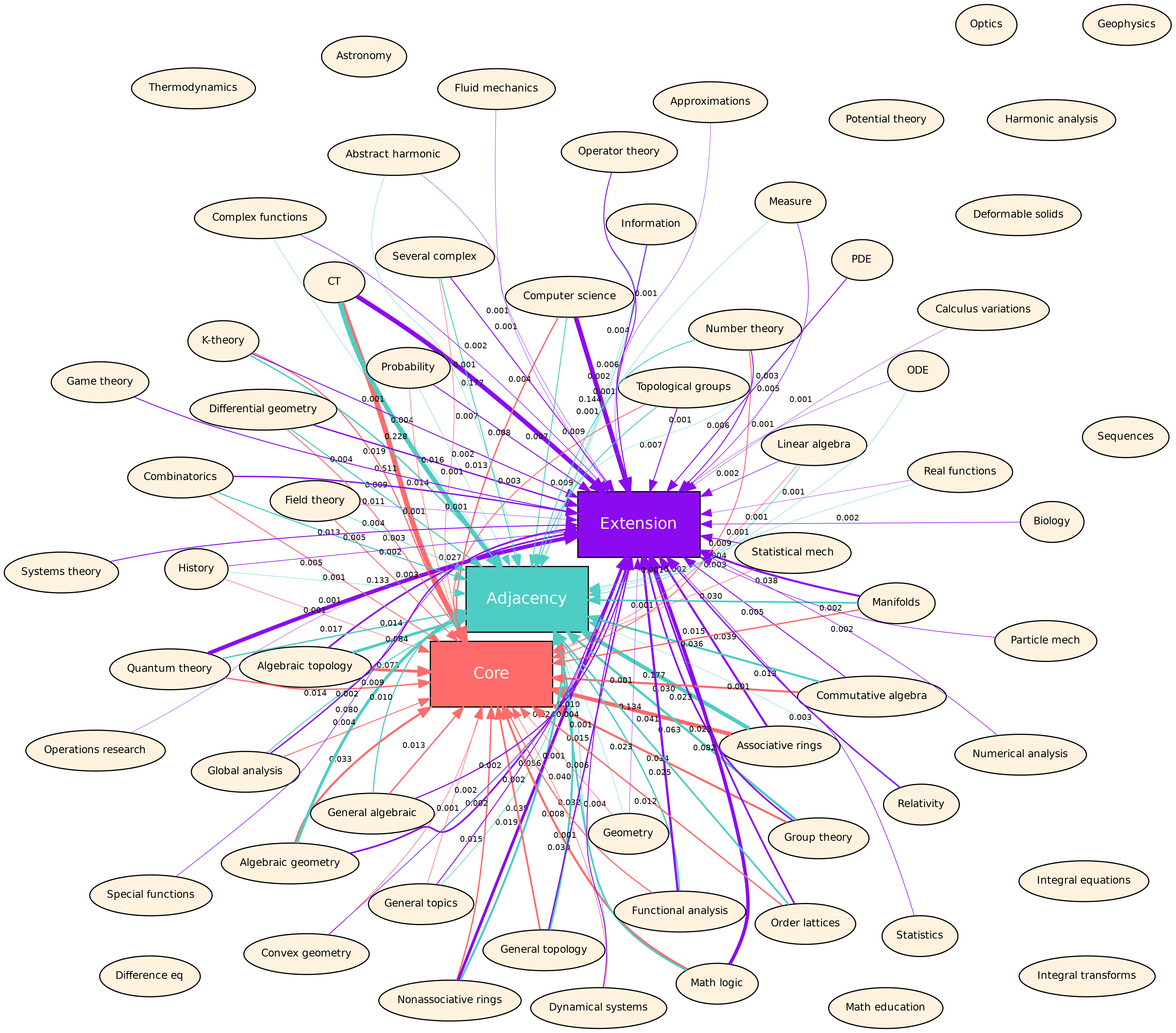}
    \caption{Knowledge flow through citations across the core, adjacency, and extension sets. Publications are classified as \emph{core} when all their MSC descriptors belong to Category Theory (MSC 18), \emph{adjacency} when their descriptors include Category Theory and Other Topics but no Applied Topics, and \emph{extension} when at least one descriptor belongs to Applied Topics (MSC 60--97). For each set, the network connects the set to the topics (three-digit MSC descriptors) associated with publications cited by documents in that set. Only citations for which MSC descriptors are available for both the citing and cited publications are included. Edge weights represent the fraction of the total frequency of each cited topic associated with each of the three sets.}
    \label{fig:citation-flow}
\end{figure}

The distinction between core, adjacency, and extension characterizes publications according to their MSC classifications, but it does not imply that these regions are informationally isolated. We therefore turn to citation patterns to examine how knowledge flows across their boundaries and, in particular, whether publications in the category-theoretic core remain connected to research outside MSC 18. Figure~\ref{fig:citation-flow} shows the network between the three sets and the MSC descriptors associated with their citations. The core is connected to 31 descriptors, compared with 40 for adjacency and 51 for extension.

Despite being defined by publications classified exclusively under Category Theory, the core exhibits substantial citation connections beyond MSC 18. Approximately 51\% of the topics associated with publications cited by the core belong to Category Theory, while the remaining 49\% are distributed across other mathematical and applied areas. Among the most prominent mathematical connections are associative rings and algebras, algebraic topology, mathematical logic, group theory, and commutative algebra.

Out of all MSC descriptors associated with publications cited by the core, 1.4\% correspond to Quantum theory and 1.3\% to computer science. Although these shares are comparatively small, their presence shows that even publications classified exclusively under Category Theory remain connected through citations to applied domains. The category-theoretic core is therefore not isolated from the broader disciplinary landscape identified in the preceding analyses: its boundaries in terms of MSC classification remain permeable at the level of citations.

Taken together, these results reveal three complementary aspects of the development of category theory. First, the field has undergone sustained growth accompanied by a substantial diversification of its topical scope. Second, this expansion has been accompanied by the growth of its research community and an increasing prevalence of collaboration, together with substantial variation in the disciplinary breadth of individual researchers. Third, interactions with neighboring mathematical and applied domains are organized through distinct categorical interfaces rather than a single generic core: theoretical computer science, quantum theory, algebra, topology, and logic connect to category theory through different combinations of categorical topics. The publication- and citation-level analyses further show that these disciplinary connections extend beyond MSC classification, with even the category-theoretic core remaining connected through citations to a broad range of mathematical and applied areas.

\section{Preliminary Conclusions}
\label{sec:conclusions}

Our analysis provides a large-scale empirical account of the development of category theory through three complementary dimensions: its historical growth and thematic diversification, the structure and communication patterns of its research community, and its interactions with neighboring mathematical and applied domains. Across these dimensions, the results portray a field that has expanded substantially while becoming increasingly diversified and connected to a broad disciplinary landscape.

\paragraph{Growth and thematic diversification.}

Category theory has undergone sustained expansion over the past six decades. Annual publication output now consistently exceeds 800 papers, while the number of active researchers has grown from a few dozen before 1960 to more than 1,500 per year. This growth has been accompanied by a substantial expansion of the field's topical scope. New MSC topics continue to appear in category theory publications, and the cumulative number of associated topics shows no clear sign of saturation. At the same time, the concentration of publications among the most frequent topics has declined. The historical development of category theory therefore reflects not simply increasing publication activity, but a progressive diversification of the mathematical subjects with which categorical research is associated.

\paragraph{A growing and increasingly collaborative research community.}

The expansion of the literature is mirrored by the development of its research community, comprising 13,392 authors in our corpus. Authorship is highly heterogeneous: nearly half of all authors are represented by one publication, while around 8\% have ten or more. Collaboration has also become increasingly prevalent. The proportion of multi-author publications has risen substantially over time, around 75\% of researchers have co-authored at least one paper, and authors with larger publication records tend to have more co-authored work.

Researchers also vary considerably in the disciplinary breadth of their publications. While most are associated with a relatively small number of major MSC areas, a small subset contributes to publications spanning ten or more mathematical and applied areas. In this sense, these researchers may act as \emph{propagators} of categorical ideas across disciplinary boundaries. The author-level results therefore reveal a community combining substantial specialization with forms of collaboration and interdisciplinary activity through which categorical research extends across different domains.

\paragraph{Distinct disciplinary interfaces and a permeable core.}

The disciplinary interactions of category theory are not organized through a single generic categorical interface. Among applied domains, theoretical computer science constitutes by far the strongest connection with category theory, followed by quantum theory, while systems and control form a smaller but clearly identifiable interface. Additional, more limited connections appear in areas including artificial intelligence, economics, network science, and biology.

These applied domains, however, do not connect to category theory in the same way. Instead, different areas interact with the field through different combinations of categorical topics. Theoretical computer science interacts particularly strongly with \emph{Categories and theories}, \emph{Special categories}, and \emph{General theory of categories and functors}. Quantum theory, by contrast, connects primarily through \emph{Categorical structures}, \emph{Monoidal categories and operads}, and \emph{Homological algebra in category theory}. A similar differentiation appears among neighboring areas of mathematics: algebra, topology, and logic each interact with category theory through distinct categorical topics. These results suggest that the circulation of categorical ideas across disciplinary boundaries is organized through specialized conceptual gateways rather than through a single categorical core.

The publication-level analysis provides a complementary view of this disciplinary structure. Nearly three quarters of publications belong to the adjacency set, combining Category Theory with other mathematical topics, while 14.2\% lie entirely within Category Theory and 13.1\% include Applied Topics. Moreover, classification within the category-theoretic core does not imply informational isolation. Citation patterns show that publications classified exclusively under Category Theory remain connected to a broad range of mathematical and applied topics. In this precise sense, the category-theoretic core is \emph{permeable}: its boundaries as defined by MSC classification do not constitute boundaries to the circulation of information through citations.

These findings provide an empirical counterpart to the historical picture with which we began. Category theory originated as a language for relating mathematical constructions across domains, and its contemporary bibliographic record continues to bear traces of this relational character. Its development is marked simultaneously by internal diversification, the expansion and increasing collaboration of its research community, and sustained interaction with neighboring mathematical and applied disciplines. Importantly, these interactions are structured rather than uniform: different domains encounter category theory through different categorical interfaces.

\subsection{Limitations and Future Directions}
\label{subsec:limitations}
These results should be interpreted in light of the limitations of the underlying data. Our analysis relies on zbMATH Open, which provides a high-quality, expert-curated view of the mathematical literature but does not fully capture the presence of category theory in adjacent disciplines. In particular, categorical research in theoretical computer science, physics, and other applied domains may appear in venues that are only partially indexed by zbMATH or may be classified in ways that do not include MSC 18. The external reach of category theory is therefore likely to be underrepresented in our corpus. More generally, bibliographic metadata capture patterns of publication, classification, authorship, and citation, but these patterns need not coincide directly with the deeper conceptual influence or transmission of mathematical ideas.

MSC classifications are periodically revised, and records may be reclassified accordingly. Consequently, some temporal patterns may reflect changes in classification practice rather than changes in underlying research activity. This effect is likely less pronounced at the coarse-grained level used throughout most of our analysis.

These limitations also suggest natural directions for extending the present analysis. Integrating additional corpora, including computer science repositories and domain-specific archives, could provide a more complete picture of the reach of categorical methods beyond mathematics. Citation directionality and temporal information could also be examined in greater detail to distinguish different forms of interaction and to investigate more directly how knowledge moves between category theory and neighboring domains.

More fundamentally, the present analysis operates primarily at the level of bibliographic entities and their relations: publications, authors, MSC topics, co-occurrences, and citations. It does not yet represent the internal conceptual structure of the mathematics itself. This distinction points toward a natural continuation of the work. The recurring topic combinations and disciplinary interfaces identified here can guide more refined analyses of the concepts, definitions, and mathematical structures through which these connections are realized. In ongoing work on what we term \emph{network mathematics}, we aim to complement these bibliographic patterns with semantic and relational representations of mathematical knowledge, moving from the traces left by a mathematical community toward a more direct representation of the ideas that circulate within it.

In this sense, scientometric analysis can serve as a guide for more refined, concept-level investigations. The identification of recurring combinations of topics, or of authors and communities that systematically bridge distinct areas, may help direct attention toward previously unnoticed correspondences. In ongoing work on network mathematics, we aim to complement these quantitative observations with semantic and relational models capable of representing mathematical concepts and their interactions more directly. The long-term goal is not only to describe the evolution of the field, but also to support the discovery of new connections within it.

\bibliographystyle{eptcs}
\bibliography{generic}

\section{Supplementary material}
\label{supplement}
\paragraph{Interactions within Other Topics}

\begin{figure}[!h]
    \centering
    \includegraphics[width=\linewidth]{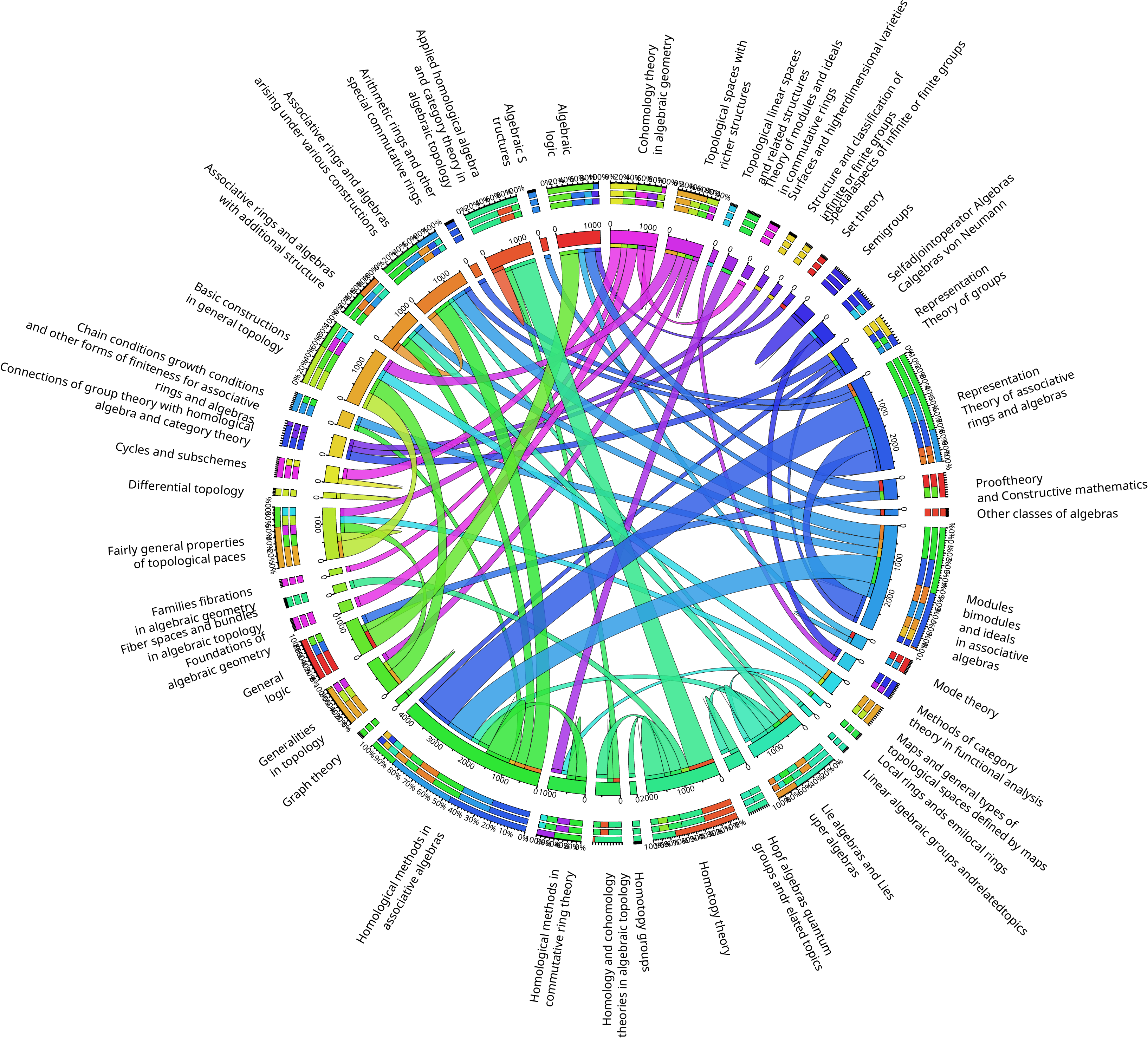}
    \caption{Chord diagram of co-occurrences between Other Topics in the zbMATH Open corpus of 33,947 Category Theory publications (MSC 18, primary or secondary classification, 1942--2024). Nodes correspond to three-digit MSC descriptors in the Other Topics class (MSC codes below 60, excluding MSC 18), and each chord represents the number of articles simultaneously classified under the corresponding pair of descriptors. The thickness of each chord is proportional to the frequency of co-occurrence. Because MSC classifications are aggregated to the three-digit level, distinct more specific descriptors belonging to the same three-digit topic may produce self-loops.}
    \label{fig:others-others}
\end{figure}

In contrast to the Applied Topics network, whose interactions are organized around the computer science hub formed by \emph{Theory of computing}, the interactions between Other Topics are distributed more evenly among the individual MSC topics, see Figure~\ref{fig:others-others} in Section \ref{supplement}. For the purposes of discussion, the strongest interactions can be grouped into a handful of broader mathematical areas. The most prominent are associative algebra, algebraic topology, general topology, logic, and algebraic geometry.

The largest of these areas is associative algebra, organized around the interaction between \emph{Homological methods in associative algebras} and \emph{Representation theory of associative rings and algebras} (959 co-occurrences), followed by interactions with \emph{Modules, bimodules and ideals in associative algebras} (791), \emph{Associative rings and algebras arising under various constructions} (487), and \emph{Associative rings and algebras with additional structure} (300). These interaction patterns place homological methods at the center of the categorical study of associative algebra. 

Algebraic topology forms another major area, anchored by the interaction between \emph{Homotopy theory} and \emph{Applied homological algebra and category theory in algebraic topology} (750), followed by \emph{Homology and cohomology theories} (318) and \emph{Fiber spaces and bundles} (229). These interactions recover the well-established interplay between category theory, homotopy theory, and homological algebra.

General topology forms another prominent area, dominated by the interactions of \emph{Basic constructions in general topology} with \emph{Fairly general properties of topological spaces} (485), \emph{Generalities in topology} (479), and \emph{Maps and general types of topological spaces} (303), together with the interaction between \emph{Generalities in topology} and \emph{Fairly general properties of topological spaces} (285).

The remaining prominent areas are logic and algebraic geometry. Logic is centered on \emph{General logic}, whose strongest interactions are with \emph{Algebraic logic} (537) and \emph{Proof theory and constructive mathematics} (258). Algebraic geometry is organized around \emph{Cohomology theory in algebraic geometry}, which interacts primarily with \emph{Cycles and subschemes} (309), \emph{Foundations of algebraic geometry} (298), and \emph{Surfaces and higher-dimensional varieties} (234).

\paragraph{Interactions of Applied Topics and Other Topics}

\begin{figure}
    \centering
    \includegraphics[width=\linewidth]{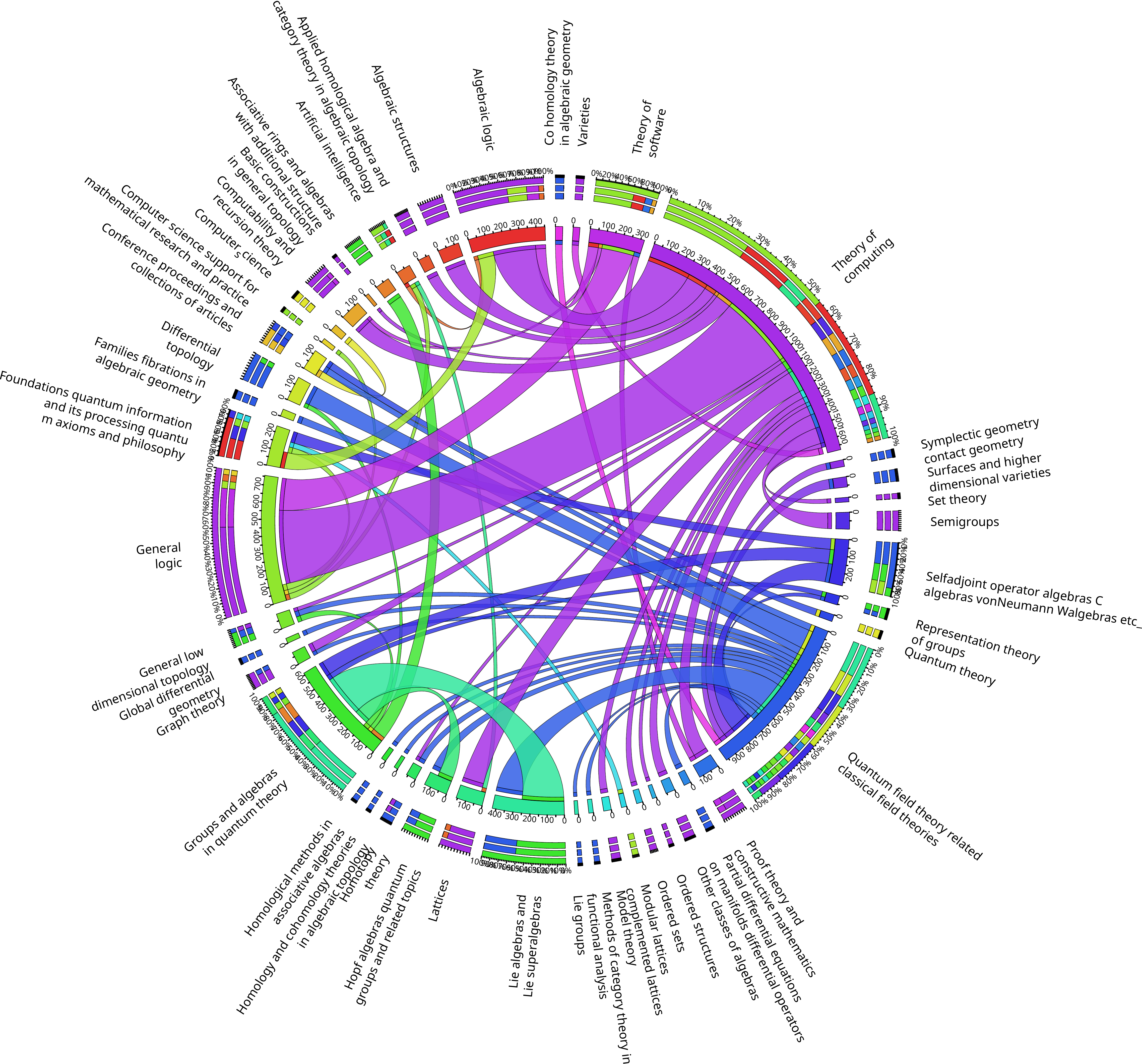}
    \caption{Chord diagram of co-occurrences between Applied Topics and Other Topics in the zbMATH Open corpus of 33,947 Category Theory publications (MSC 18, primary or secondary classification, 1942--2024). Nodes correspond to three-digit MSC descriptors in Applied Topics (MSC codes 60--97) and Other Topics (MSC codes below 60, excluding MSC 18), and each chord represents the number of articles simultaneously classified under the corresponding pair of descriptors, with one descriptor from each class. The thickness of each chord is proportional to the frequency of co-occurrence.}
    \label{fig:app-others}
\end{figure}

This section complements the previous analyses by showing the direct interactions between Applied Topics and Other Topics within the category theory corpus (see Figure~\ref{fig:app-others} in Section \ref{supplement}). It reveals which branches of mathematics are most closely associated with different application domains, exposing the mathematical foundations that support the applications represented in the corpus. 

The interactions between Applied Topics and Other Topics are dominated by two major interfaces: theoretical computer science with logic, and quantum theory with algebra and geometry. 

The strongest interaction in the entire diagram is between \emph{Theory of computing} and \emph{General logic} (458 co-occurrences), followed by \emph{Algebraic logic} (265), \emph{Lattices} (132), \emph{Algebraic structures} (126), \emph{Semigroups} (99), \emph{Computability and recursion theory} (98), and \emph{Proof theory and constructive mathematics} (87). These interactions identify logic and foundational algebraic structures as the principal mathematical interface of theoretical computer science.

The second major interface is quantum theory. \emph{Groups and algebras in quantum theory} interacts most strongly with \emph{Lie algebras and Lie superalgebras} (256), followed by \emph{Hopf algebras, quantum groups and related topics} (90), \emph{Associative rings and algebras with additional structure} (85), and \emph{Selfadjoint operator algebras} (88). Likewise, \emph{Quantum field theory and related classical field theories} is connected primarily with \emph{Lie algebras and Lie superalgebras} (166), \emph{Differential topology} (118), \emph{Selfadjoint operator algebras} (115), and \emph{Surfaces and higher-dimensional varieties} (60).

Several smaller interfaces are also evident. \emph{Theory of computing} interacts with \emph{Applied homological algebra and category theory in algebraic topology} (67), \emph{Graph theory} (53), and \emph{Model theory} (45), while \emph{Theory of software} also exhibits substantial interactions with \emph{General logic} (200) and \emph{Algebraic logic} (63).

\end{document}